\documentclass[12pt]{amsart}

\usepackage{amsmath,amssymb,amscd}

\input{prepictex}
\input{pictex}
\input{postpictex}

\def\arr(#1,#2)(#3,#4){\arrow <0.15cm> [0.25, 0.75] from {#1} {#2} to {#3} {#4}}

\def\NN{\N}%{\mathbf N}
\def\TT{\T}%{\mathbf T}
\def\ZZ{\Z}%{\mathbf Z}
\def\Mm{\mathcal M}

\def\ker{\operatorname{ker}}

\def\Obj{\operatorname{Obj}}
\def\Mor{\operatorname{Mor}}
\def\id{\operatorname{id}}

\def\Aut{\operatorname{Aut}}
\def\lsp{\operatorname{span}}
\def\clsp{\operatorname{\overline{span\!}\,\,}}
\def\mod{\operatorname{mod}}

\theoremstyle{plain}
\newtheorem{thm}{Theorem}[section]
\newtheorem*{thm*}{Theorem}
\newtheorem{prop}[thm]{Proposition}
\newtheorem*{prop*}{Proposition}
\newtheorem{lemma}[thm]{Lemma}
\newtheorem{cor}[thm]{Corollary}
\newtheorem*{conj*}{Conjecture}
\newtheorem*{cor*}{Corollary}
\newtheorem{defn}[thm]{Definition}

\theoremstyle{definition}
\newtheorem*{defn*}{Definition}
\newtheorem{rems}[thm]{Remarks}
\newtheorem*{rems*}{Remarks}
\newtheorem*{proof*}{Proof}
\newtheorem*{not*}{Notation}
\newtheorem{notation}[thm]{Notation}

 \newcommand{\A}{{\mathcal A}}
        \newcommand{\D}{{\mathcal D}}\newcommand{\HH}{{\mathcal H}}
        \newcommand{\LL}{{\mathcal L}}
        \newcommand{\B}{{\mathcal B}}
        \newcommand{\K}{{\mathcal K}}

        \newcommand{\s}{\sigma}
\newcommand{\al}{\alpha}
        \newcommand{\p}{\partial}
        \newcommand{\dd}{|\D|}
\newcommand{\bma}{\left(\begin{array}{cc}}
\newcommand{\ema}{\end{array}\right)}
\newcommand{\bca}{\left(\begin{array}{c}}
\newcommand{\eca}{\end{array}\right)}
\def\clsp{\overline{\operatorname{span}}}
\def\T{\mathbf T}
\def\Aut{\operatorname{Aut}}

        \newcommand{\Q}{\mathbf Q}
        \newcommand{\R}{\mathbf R}
        \newcommand{\C}{\mathbf C}
        
\newcommand{\Z}{\mathbf Z}
\newcommand{\N}{\mathbf N}

\newcommand{\ben}{\begin{displaymath}}
        \newcommand{\een}{\end{displaymath}}
\newcommand{\be}{\begin{equation}}
\newcommand{\ee}{\end{equation}}

        \newcommand{\bean}{\begin{eqnarray*}}
        \newcommand{\eean}{\end{eqnarray*}}
\newcommand{\nno}{\nonumber\\}
\newcommand{\bea}{\begin{eqnarray}}
        \newcommand{\eea}{\end{eqnarray}}

\def\cross#1{\rlap{\hskip#1pt\hbox{$-$}}}
        
        \def\bigintcross{\cross{2.3}\int}

\newcommand{\nc}{\newcommand}
\nc{\nt}{\newtheorem} \nc{\gf}[2]{\genfrac{}{}{0pt}{}{#1}{#2}}
\nc{\mb}[1]{{\mbox{$ #1 $}}} \nc{\real}{{\mathbf R}}
\nc{\comp}{{\mathbf C}} \nc{\ints}{{\mathbf Z}}
\nc{\Ltoo}{\mb{L^2({\mathbf H})}} \nc{\rtoo}{\mb{{\mathbf R}^2}}
\nc{\slr}{{\mathbf {SL}}(2,\real)} \nc{\slz}{{\mathbf
{SL}}(2,\ints)} \nc{\su}{{\mathbf {SU}}(1,1)} \nc{\so}{{\mathbf
{SO}}} \nc{\hyp}{{\mathbb H}} \nc{\disc}{{\mathbf D}}
\nc{\torus}{{\mathbf T}}

\nc{\ca}{{\mathcal A}} \nc{\cag}{{{\mathcal A}^\Gamma}}
\nc{\cg}{{\mathcal G}} \nc{\chh}{{\mathcal H}} \nc{\ck}{{\mathcal
B}} \nc{\cl}{{\mathcal L}} \nc{\cm}{{\mathcal M}}
\nc{\cn}{{\mathcal N}} \nc{\cs}{{\mathcal S}} \nc{\cz}{{\mathcal Z}}
\nc{\sind}{\sigma{\rm -ind}}
\newcommand{\la}{\langle}
\newcommand{\ra}{\rangle}
\def\Ends{\operatorname{Ends}}
\title{The noncommutative geometry of $k$-graph $C^*$-algebras}
\author{David Pask}
\author{Adam Rennie}
\author{Aidan Sims}
\address{Department of Mathematics\\
University of Newcastle\\
NSW  2308\\
AUSTRALIA}\address{Institute for Mathematical Sciences\\Universitetsparken
5\\DK-2100\\Copenhagen\\Denmark}
\email{david.pask@newcastle.edu.au}
\email{rennie@math.ku.dk}
\email{aidan.sims@newcastle.edu.au} 
\date{December 23, 2005}
\thanks{The last two authors were supported by Australian
Research Council Australian Postdoctoral Fellowships. In addition, AR
was supported by Statens Naturvidenskabelige Forskningsr{\aa}d, Denmark}
\begin{document}

\begin{abstract}
This paper is comprised of two related parts.
First we discuss which $k$-graph algebras have faithful
traces. We characterise the existence of a faithful semifinite
lower-semicontinuous gauge-invariant trace on $C^*(\Lambda)$ in terms
of the existence of a faithful graph trace on $\Lambda$.

Second, for $k$-graphs with faithful gauge invariant trace, we construct a
smooth $(k,\infty)$-summable semifinite spectral triple. We use
the semifinite local index theorem to compute the pairing with
$K$-theory. This numerical pairing can be obtained by applying the
trace to a $KK$-pairing with values in the $K$-theory of the fixed
point algebra of the $\T^k$ action. As with graph algebras, the index pairing is an
invariant for a finer structure than the isomorphism class of the algebra.
\end{abstract}

\maketitle

\section{Introduction}

In this paper we generalise the construction of semifinite spectral
triples for graph algebras of \cite{PRen} to the $C^*$-algebras of
higher-rank graphs, or $k$-graphs. Experience with $k$-graph
algebras has shown that from a $C^*$-algebraic point of view they
tend to behave very much like graph $C^*$-algebras. Consequently the
transition from graph $C^*$-algebras to $k$-graph $C^*$-algebras
often appears quite simple. The subtlety generally lies in the added
combinatorial complexity of $k$-graphs, and in particular in
identifying the right higher-dimensional analogues of the
graph-theoretic conditions which arise in the one-dimensional case.
This experience is borne out again in the current paper: once the
appropriate $k$-graph theoretic conditions have been identified, the
generalisations of the constructions in \cite{PRen} to higher-rank
graphs turn out to be mostly straightforward.

However, the pay-offs from carrying through this analysis are significant 
from the points of view of both noncommutative geometry and $k$-graph 
algebra theory.

The pay-off for noncommutative geometry is that our analysis allows
us to construct infinitely many examples of (semifinite) spectral
triples of every integer dimension $k\geq 1$ ($k=1$ is contained in
\cite{PRen}). These spectral triples are generically semifinite, and
so come from $KK$ classes rather than $K$-homology classes.
Computations can be made very explicitly with these algebras, and we
use this to relate the semifinite index pairing to the $KK$-index
for these examples. This has led to a general picture of the
relationship between semifinite index theory and $KK$-theory,
\cite{KNR}. From our point of view, $k$-graph algebras are
sufficiently generic to reveal the relationship between $KK$-theory
and semifinite index theory. We highlight this at the end of
Section~\ref{index}.

The connections between non-commutative geometry and classical
(commutative) differential geometry are still a subject of intense
research (and speculation). This paper helps to illuminate these 
connections in two ways. The first, and always the most important for
understanding the extension of geometry to the noncommutative
setting, is index theory. We relate an analytic index to a
$KK$-index with values in the $K$-theory of the fixed point algebra.
This $KK$-index is far better adapted to topological interpretation, and
furthermore exists in greater generality.

The second connection between this paper and classical geometry is less
deep, but is a key tool in our construction. The group $\T^k$ acting
on a $k$-graph algebra is a Lie group. We use the action of $\T^k$
to `push forward' the Dirac operator on $\T^k$ to a Dirac operator
for a $k$-graph algebra. While we have done this only for the (spin)
Dirac operator of the simplest spin structure on $\T^k$, the
possibility also exists for repeating our construction for every
$\operatorname{spin}^c$ structure on $\T^k$, as well as the Hodge-de
Rham operator, and other twisted Dirac operators. It is an
interesting question to determine to what extent this `pushing
forward' operation can be systematised.

The pay-off for $k$-graph algebra theory is that we obtain
detailed information about semifinite traces on $k$-graph algebras.
To construct a semifinite spectral triple one requires a faithful 
semifinite trace on the underlying $C^*$-algebra. Hence our first
step is to investigate when such a trace exists on a $k$-graph algebra.
We characterise faithful semifinite gauge-invariant traces on 
$C^*(\Lambda)$ in terms of \emph{graph traces} on $\Lambda$ 
(cf. \cite{H, PRen}). In particular, this
represents a first systematic exploration of gauge-invariant traces
on $C^*$-algebras associated to graphs of arbitrary rank (see
Section~\ref{traces}).

In the appendix, we also identify a (fairly restrictive) class of
$k$-graphs which admit faithful graph traces, demonstrate that their
$C^*$-algebras are Morita equivalent to direct sums of algebras of
continuous functions on tori of rank $0, \dots, k$, and calculate
their $K$-theory. For graph algebras this is not new because the
$K$-theory of graph algebras is completely understood \cite{RSz}.
However, only for $2$-graphs have general $K$-theory computations
recently emerged \cite{EPhD}. Consequently any advances on
$K$-theory for general $k$-graph algebras are significant.

\textbf{Outline.} The paper is arranged a follows. In
Section~\ref{kgraph} we review the basic definitions of $k$-graphs
and $k$-graph algebras. In Section~\ref{traces} we show that a 
$k$-graph algebra $C^*(\Lambda)$ admits a faithful, semifinite, 
lower-semicontinuous, gauge-invariant trace if and only if $\Lambda$
admits a faithful graph trace (in the Appendix, we 
identify a substantial class of $k$-graphs which admit such a graph 
trace, show that they are Morita equivalent to commutative 
$C^*$-algebras, and hence compute their $K$-theory). 

Section~\ref{spectrip} reviews the definitions we require pertaining
to semifinite spectral triples. In Section~\ref{firstconstruction}
we construct a Kasparov module for the $C^*$-algebra of any locally
finite, locally convex $k$-graph with no sinks. This is a very
general construction, and the resulting Kasparov module is even iff
$k$ is an even integer. In Section \ref{secondconstruction} we
construct $(k,\infty)$-summable spectral triples for $k$-graph
algebras with faithful trace.

The von Neumann algebra constructed as part of our spectral triple plays
the role of the crossed product of the (von Neumann completion of the)
graph algebra by the $\T^k$ action. The precise relationship is not clear to
us, but we suspect the two are isomorphic.

In Section~\ref{index} we use these spectral triples to compute
index pairings and compare them with the Kasparov product as in
\cite{PRen}. By example we indicate how the semifinite index can be
used to obtain more refined information than the usual Fredholm
index.

\textbf{Acknowledgements} The second author would like to thank Ian
Putnam for bringing the $K$-theory constructions of Section
\ref{index} to his attention. We would also like to thank our
referee for a careful reading of the manuscript and detailed
comments leading to significant improvements in the exposition.

\section{$k$-Graph $C^*$-Algebras}\label{kgraph}

\subsection{Higher-rank graphs and their $C^*$-algebras}
In this subsection we outline the basic notation and definitions
of $k$-graphs and their $C^*$-algebras. We refer the reader to
\cite{RSY1} for a more thorough account.

\vskip3pt \textbf{Higher-rank graphs.} Throughout this paper,
we regard $\N^k$  as a monoid under pointwise addition. We
denote the usual generators of $\N^k$ by $e_1, \dots, e_k$, and
for $n \in \N^k$ and $1 \le i \le k$, we denote the $i^{\rm th}$
coordinate of $n$ by $n_i \in \N$; so $n = \sum n_i e_i$.
For $m,n \in \N^k$, we write $m \le n$ if $m_i \le n_i$ for all
$i$. By $m < n$, we mean $m \le n$ and $m \not= n$. We use $m \vee
n$ and $m \wedge n$ to denote, respectively, the coordinate-wise
maximum and coordinate-wise minimum of $m$ and $n$; so that $m
\wedge n \le m, n \le m \vee n$ and these are respectively the
greatest lower bound and least upper bound of $m,n$ in $\N^k$.

\begin{defn}[Kumjian-Pask~\cite{KP}]
A \emph{graph of rank $k$} or \emph{$k$-graph} is a pair
$(\Lambda,d)$ consisting of a countable category $\Lambda$ and a
\emph{degree} functor $d : \Lambda \to \N^k$ which satisfy the
following \emph{factorisation property}: if $\lambda \in
\Mor(\Lambda)$ satisfies $d(\lambda) = m+n$, then there are unique
morphisms $\mu,\nu \in \Mor(\Lambda)$ such that $d(\mu) = m$,
$d(\nu) = n$, and $\lambda = \mu\circ\nu$.
\end{defn}

The factorisation property ensures (see \cite{KP}) that the
identity morphisms of $\Lambda$ are precisely the morphisms of
degree 0; that is $\{\id_o : o \in \Obj(\Lambda)\} = d^{-1}(0)$.
This means that we may identify each object with its identity
morphism, and we do this henceforth. This done, we can regard
$\Lambda$ as consisting only of its morphisms, and we write
$\lambda \in \Lambda$ to mean $\lambda \in \Mor(\Lambda)$.

Since we are thinking of $\Lambda$ as a kind of graph, we write
$r$ and $s$ for the codomain and domain maps of $\Lambda$
respectively. We refer to elements of $\Lambda$ as \emph{paths},
and to the paths of degree $0$ (which correspond to the objects of
$\Lambda$ as above) as \emph{vertices}. Extending these
conventions, we refer to the elements of $\Lambda$ with minimal
nonzero degree (that is $d^{-1}(\{e_1,\dots, e_k\})$ as
\emph{edges}.

\begin{notation}
To try to minimise confusion, we will always use $u,v,w$ to denote
vertices, $e,f$ to denote edges, and lower-case Greek letters
$\lambda,\mu,\nu$, etc{.} for arbitrary paths. We will also drop
the composition symbol, and simply write $\mu\nu$ for $\mu \circ
\nu$ when the two are composable.
\end{notation}
\textsc{Warning:} \textsl{because $\Lambda$ is a category, composition of
morphisms reads from right to left. Hence paths $\mu$ and $\nu$ in
$\Lambda$ can be composed to form $\mu\nu$ if and only if $r(\nu)
= s(\mu)$, and in this case, $r(\mu\nu) = r(\mu)$ and $s(\mu\nu) =
s(\nu)$. This is the reverse of the  convention for
directed graphs, used in \cite{BPRS,kpr,KPRR,PRen},
so the reader should beware. In particular the
roles of sources and sinks, and of ranges and sources,
 are opposite to those in \cite{PRen}.}
\begin{defn} For each $n \in \N^k$, we write $\Lambda^n$ for the collection
$\{\lambda \in \Lambda : d(\lambda) = n\}$ of paths of degree $n$.
\end{defn}
The range and source  $r,s$ are thus maps from $\Lambda$
to $\Lambda^0$, and if $v \in \Lambda^0$, then $r(v) = v = s(v)$.

Given $\lambda \in \Lambda$ and $S \subset \Lambda$, it makes
sense to write $\lambda S$ for $\{\lambda\sigma : \sigma \in S,
r(\sigma) = s(\lambda)\}$, and likewise $S\lambda =
\{\sigma\lambda : \sigma \in S, s(\sigma) = r(\lambda)\}$. In
particular, if $v \in \Lambda^0$, then $v S$ is the collection of
elements of $S$ with range $v$, and $S v$ is the collection of
elements of $S$ with source $v$.
\begin{defn}
Let $(\Lambda,d)$ be a $k$-graph. We say that $\Lambda$ is
\emph{row-finite} if $|v\Lambda^n| < \infty$ for each $v \in
\Lambda^0$ and $n \in \N^k$. We say that $\Lambda$ is
\emph{locally-finite} if it is row-finite and also satisfies
$|\Lambda^n v| < \infty$ for all $v \in \Lambda^0$ and $n \in
\N^k$. We say that $\Lambda$ has \emph{no sources} (resp.
\emph{no sinks}) if $v\Lambda^n$ (resp. $\Lambda^n v$) is nonempty
for each $v \in \Lambda^0$ and $n \in \N^k$. Finally, we say that
$\Lambda$ is \emph{locally convex} if, for each edge $e \in
\Lambda^{e_i}$, and each $j \not= i$ in $\{1,\dots,k\}$, we have
$s(e)\Lambda^{e_j} = \emptyset$ only if $r(e)\Lambda^{e_j} =
\emptyset$.
\end{defn}

As in \cite{RSY1}, for locally convex $k$-graphs, we
use the notation $\Lambda^{\le n}$ to denote the collection
\[
\Lambda^{\leq n}:=\{\lambda \in \Lambda : d(\lambda) \le n, \mu \in
s(\lambda)\Lambda\text{ and } d(\lambda\mu) \le n \text{ implies }
\mu = s(\lambda)\}.
\]
Intuitively, $\Lambda^{\le n}$ is the collection of paths whose
degree is ``as large as possible" subject to being dominated by
$n$. In a $1$-graph, $\Lambda^{\le n}$ is the set of paths
$\lambda \in \Lambda$ whose length is at most $n$ and is less than
$n$ only if $s(\lambda)$ receives no edges. The significance of
this is that the partial isometries associated to distinct paths in
$\Lambda^{\le n}$ have orthogonal range projections (cf. relation~(CK4)
below). For more on the importance of $\Lambda^{\le n}$, see \cite{RSY1}.

\begin{figure}[ht]
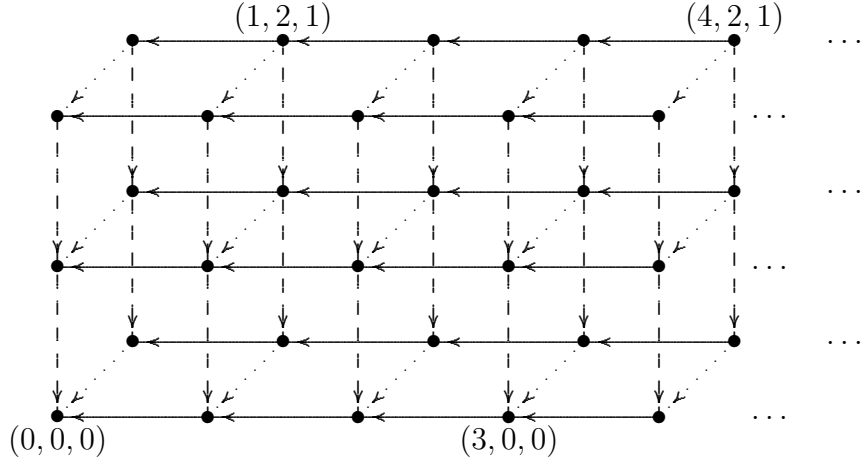

\[{
\def\xArr{
\arr(.3,0)(.2,0) }
\def\yArr{
\arr(0,.3)(0,.2) }
\def\zArr{
\arr(.2,.2)(.199,.199) }
\beginpicture
\setcoordinatesystem units <1cm, 1cm> \put{\put{$\bullet$} at 0 0
\xArr\yArr\zArr}[lB] at 0 0 \put{$(0,0,0)$}[ct] at 0 -.1
\put{\put{$\bullet$} at 0 0 \xArr\yArr\zArr}[lB] at 2 0
\put{\put{$\bullet$} at 0 0 \xArr\yArr\zArr}[lB] at 4 0
\put{\put{$\bullet$} at 0 0 \xArr\yArr\zArr}[lB] at 6 0
\put{$(3,0,0)$}[ct] at 6 -.1 \put{\put{$\bullet$} at 0 0
\yArr\zArr}[lB] at 8 0 \put{$\dots$} at 9.5 0 \plot 0 0 8 0 /
\put{\put{$\bullet$} at 0 0 \xArr\yArr\zArr}[lB] at 0 2
\put{\put{$\bullet$} at 0 0 \xArr\yArr\zArr}[lB] at 2 2
\put{\put{$\bullet$} at 0 0 \xArr\yArr\zArr}[lB] at 4 2
\put{\put{$\bullet$} at 0 0 \xArr\yArr\zArr}[lB] at 6 2
\put{\put{$\bullet$} at 0 0 \yArr\zArr}[lB] at 8 2 \put{$\dots$}
at 9.5 2 \plot 0 2 8 2 /
\put{\put{$\bullet$} at 0 0 \xArr\zArr}[lB] at 0 4
\put{\put{$\bullet$} at 0 0 \xArr\zArr}[lB] at 2 4
\put{\put{$\bullet$} at 0 0 \xArr\zArr}[lB] at 4 4
\put{\put{$\bullet$} at 0 0 \xArr\zArr}[lB] at 6 4
\put{\put{$\bullet$} at 0 0 \zArr}[lB] at 8 4 \put{$\dots$} at 9.5
4 \plot 0 4 8 4 / \setdashes \plot 0 4 0 0 / \plot 2 4 2 0 / \plot
4 4 4 0 / \plot 6 4 6 0 / \plot 8 4 8 0 / \setsolid
%
% Second array
\put{\put{$\bullet$} at 0 0 \xArr\yArr}[lB] at 1 1
\put{\put{$\bullet$} at 0 0 \xArr\yArr}[lB] at 3 1
\put{\put{$\bullet$} at 0 0 \xArr\yArr}[lB] at 5 1
\put{\put{$\bullet$} at 0 0 \xArr\yArr}[lB] at 7 1
\put{\put{$\bullet$} at 0 0 \yArr}[lB] at 9 1 \put{$\dots$} at
10.5 1 \plot 1 1 9 1 /
\put{\put{$\bullet$} at 0 0 \xArr\yArr}[lB] at 1 3
\put{\put{$\bullet$} at 0 0 \xArr\yArr}[lB] at 3 3
\put{\put{$\bullet$} at 0 0 \xArr\yArr}[lB] at 5 3
\put{\put{$\bullet$} at 0 0 \xArr\yArr}[lB] at 7 3
\put{\put{$\bullet$} at 0 0 \yArr}[lB] at 9 3 \put{$\dots$} at
10.5 3 \plot 1 3 9 3 /
\put{\put{$\bullet$} at 0 0 \xArr}[lB] at 1 5 \put{\put{$\bullet$}
at 0 0 \xArr}[lB] at 3 5 \put{$(1,2,1)$}[cb] at 3 5.1
\put{\put{$\bullet$} at 0 0 \xArr}[lB] at 5 5 \put{\put{$\bullet$}
at 0 0 \xArr}[lB] at 7 5 \put{\put{$\bullet$} at 0 0 }[lB] at 9 5
\put{$(4,2,1)$}[cb] at 9 5.1 \put{$\dots$} at 10.5 5 \plot 1 5 9 5
/ \setdashes \plot 1 5 1 1 / \plot 3 5 3 1 / \plot 5 5 5 1 / \plot
7 5 7 1 / \plot 9 5 9 1 / \setsolid
\setdots \plot 1 1 0 0 / \plot 3 1 2 0 / \plot 5 1 4 0 / \plot 7 1
6 0 / \plot 9 1 8 0 /
\plot 1 3 0 2 / \plot 3 3 2 2 / \plot 5 3 4 2 / \plot 7 3 6 2 /
\plot 9 3 8 2 /
\plot 1 5 0 4 / \plot 3 5 2 4 / \plot 5 5 4 4 / \plot 7 5 6 4 /
\plot 9 5 8 4 / \setsolid
\endpicture
}
\]
\caption{The skeleton of $\Omega_{3, (\infty,2,1)}$}
\label{fig:Omega skeleton}
\end{figure}

\vskip3pt \textbf{$\mathbf{\Omega_{k,m}}$ and boundary paths.} For
$k \ge 1$ and $m \in (\N \cup \{\infty\})^k$, we define a
$k$-graph $\Omega_{k,m}$ as follows:
\begin{gather*}
\Omega_{k,m}^0 = \{n \in \N^k : n \le m\} \qquad
\Omega_{k,m}^n = \{(p,q) \in \N^k : p,q \in \Omega_{k,m}^0, q - p = n\} \\
r(p,q) = p, \qquad s(p,q) = q, \qquad (p,q)\circ(q,n) = (p,n).
\end{gather*}
See Figure~\ref{fig:Omega skeleton} for a ``picture" of
$\Omega_{3,(\infty,2,1)}$.

Each path $\lambda$ of degree $p$ in a $k$-graph $\Lambda$ determines a
degree-preserving functor $\hat\lambda$ from $\Omega_{k,p}$ to $\Lambda$
as follows: the image $\hat\lambda(m,n)$ of the morphism
$(m,n) \in \Omega_{k,p}$ is the unique morphism in $\Lambda^{n-m}$ such
that there exist $\mu \in \Lambda^{m}$ and $\nu \in \Lambda^{p-n}$
satisfying $\lambda = \mu \hat\lambda(m,n) \nu$. (The existence
and uniqueness of $\hat\lambda(m,n)$ is guaranteed by the factorisation
property).

In fact for each $p \in \N^k$, the map $\lambda \mapsto \hat\lambda$ is
a bijection between $\Lambda^p$ and the set of degree-preserving functors
from $\Omega_{k,p}$ to $\Lambda$. In practise, we
just write $\lambda(m,n)$ for the segment $\hat\lambda(m,n)$ of $\lambda$
described in the previous paragraph, and we write $\lambda(m)$ for the
range of the path $\lambda(m,n)$, which we think of as the vertex on
$\lambda$ at position $m$. If $\lambda \in \Lambda^p$
and $0 \le m \le n \le p$, then
\[
\lambda = \lambda(0,m)\lambda(m,n)\lambda(n,p),
\quad
s(\lambda(m,n)) = \lambda(n)
\quad\text{and}\quad
r(\lambda(m,n)) = \lambda(m).
\]

We extend this correspondence between paths and degree-preserving functors
to define the notion of a boundary path in a $k$-graph.

\begin{defn} A \emph{boundary path} of a $k$-graph $\Lambda$ is a
degree-preserving functor $x : \Omega_{k,m} \to \Lambda$ such that
\[
\text{if $m_i < \infty$, $n \in \N^k$, $n \le m$ and
$n_i = m_i$, then $x(n)\Lambda^{e_i} =
\emptyset$;}
\]
so the directions in which $x$ is finite
are those in which it cannot be extended. If $x : \Omega_{k,m} \to
\Lambda$ is a boundary path,  we denote $m$ by $d(x)$, and
$x(0)$ by $r(x)$. We write $\Lambda^{\le \infty}$ for the
set of all boundary paths of $\Lambda$.
\end{defn}

Note that if $\lambda \in \Lambda$ satisfies
$s(\lambda)\Lambda^{n} = \emptyset$ for all $n > 0$ (that is, if $s(\lambda)$
is a source in $\Lambda$), then the graph morphism $\hat\lambda :
\Omega_{k, d(\lambda)} \to \Lambda$ discussed above belongs to
$\Lambda^{\le \infty}$; we think of $\lambda$ itself as a boundary path of
$\Lambda$.

\begin{defn}\label{dfn:ends}
An \emph{end} of $\Lambda$ is a boundary path $x \in \Lambda^{\le\infty}$ such
that for all $n \le d(x)$, $r(x)\Lambda^n = \{x(0,n)\}$. We denote the set of
ends of $\Lambda$ by $Ends(\Lambda)$.
\end{defn}

\begin{rems}\label{rmk:uptolambda on ends}
If $x$ is an end of $\Lambda$, then
$r(x) \Lambda^{\le n} = \{x(0, n \wedge d(x))\}$
for all $n \in \N^k$.
\end{rems}

\vskip3pt \textbf{Skeletons.} To draw a $k$-graph, we use its
\emph{skeleton}. The skeleton of a $k$-graph $\Lambda$ is the
directed graph whose vertices and edges are those of $\Lambda$,
but with the $k$ different types of edges distinguished using
$k$ different colours. In this paper, we use solid lines for edges
of degree $e_1$, dashed lines for edges of degree $e_2$, and
dotted lines for edges of degree $e_3$. For example, the skeleton
of $\Omega_{3, (\infty, 2, 1)}$ is presented in
Figure~\ref{fig:Omega skeleton}

The factorisation property says that if $e$ and $f$ are edges of
degree $e_i$ and $e_j$ respectively such that $s(e) = r(f)$, then
the path $ef$ can be expressed in the form $f'e'$ where $d(f') =
e_j$ and $d(e') = e_i$. In the skeleton for $\Omega_{3,
(\infty,2,1)}$ there is just one way this can happen; so the
skeleton is actually a commuting diagram in the category, and
although there appear to be many ways to get from $(1,2,1)$ to
$(0,0,0)$, for example, each of these paths
yields the same morphism in the category, so there is really just
one path in $\Omega_{3,
(\infty,2,1)}$ from $(1,2,1)$ to $(0,0,0)$.

The information determining the factorisation
property is not always included in the skeleton, and it must then be
specified separately as a set of \emph{factorisation
rules}. The uniqueness of factorisations ensures that amongst the
factorisation rules for the skeleton of a $k$-graph, each composition
$ef$ where $e$ and $f$ are composable edges of different colours will appear
exactly once. A set of factorisation rules for a skeleton with
this property is referred to as an \emph{allowable}
factorisation regime.

\begin{figure}[ht]
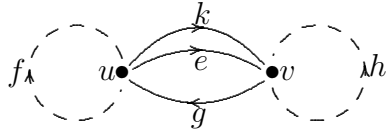

\[
\beginpicture
\setcoordinatesystem units <1cm,1cm>
\put{$\bullet$} at 4 0
\put{$\bullet$} at 6 0
\setdashes
\circulararc 320 degrees from 4 .2 center at 3.4 0
\arrow <0.15cm> [0.25,0.75] from 2.77  .01 to 2.77 .03
\circulararc 320 degrees from 6 -.2 center at 6.6 0
\arrow <0.15cm> [0.25,0.75] from 7.227 .075 to 7.226 .095
\setsolid
\setquadratic
\plot 4.15 .05 5 .3 5.85 .05 /
\arrow <0.15cm> [0.25,0.75] from 5.02 .3 to 5.06 .3
\plot 4.1 .1 5 .6 5.9 .1 /
\arrow <0.15cm> [0.25,0.75] from 5.04 .6 to 5.08 .6
\setquadratic
\plot 4.1 -.1 5 -.4 5.9 -.1 /
\arrow <0.15cm> [0.25,0.75] from 4.92 -.4 to 4.88 -.4
\put{$u$} at 3.8  0
\put{$v$} at 6.2 0
\put{$f$} at 2.6 0
\put{$e$} at 5.05 .13
\put{$k$} at 5.07 .8
\put{$g$} at 5.02 -.6
\put{$h$} at 7.4 .1
\endpicture
\]
\caption{A skeleton that admits two distinct factorisation regimes}
\label{fig:factorisations}
\end{figure}

For example, in the $1$-skeleton of Figure~\ref{fig:factorisations}
the allowable factorisation regimes are: $\{ef = he, kf = hk\}$
and $\{ef = hk, kf = he\}$).

A skeleton together with an allowable  factorisation regime
determines at most one $k$-graph. When $k = 2$, each skeleton and
allowable  factorisation regime determines a unique
$k$-graph. For $k \ge 3$, there is an additional associativity
condition on the factorisation rules which must be verified
\cite{FS1}; but the issue does not arise in the examples we give
in this paper.

\vskip3pt \textbf{Cuntz-Krieger families and $C^*(\Lambda)$.} As
with directed graphs, we are interested in higher-rank graphs
because we can associate to each one a $C^*$-algebra of
Cuntz-Krieger type.

\begin{defn}
Let $(\Lambda,d)$ be a row-finite locally convex $k$-graph. A
Cuntz-Krieger $\Lambda$-family is a collection $\{s_\lambda :
\lambda \in \Lambda\}$ of partial isometries satisfying
\begin{itemize}
\item[(CK1)] $\{s_v : v \in \Lambda^0\}$ is a collection of
mutually orthogonal projections; \item[(CK2)] $s_\mu s_\nu =
s_{\mu\nu}$ for all $\mu,\nu \in \Lambda$ with $s(\mu) = r(\nu)$;
\item[(CK3)] $s^*_\lambda s_\lambda = s_{s(\lambda)}$ for all
$\lambda \in \Lambda$; and \item[(CK4)] $s_v = \sum_{\lambda \in
v\Lambda^{\le n}} s_\lambda s^*_\lambda$ for all $v \in \Lambda^0$
and $n \in \N^k$.
\end{itemize}
As a point of notation, we will henceforth denote the vertex
projection $s_v$ by $p_v$  to remind ourselves that it is a
projection.
\end{defn}

The Cuntz-Krieger algebra of $\Lambda$, denoted $C^*(\Lambda)$, is
the universal $C^*$-algebra generated by a Cuntz-Krieger family
$\{s_\lambda : \lambda \in \Lambda\}$. By this we mean that given
any other Cuntz-Krieger $\Lambda$-family
$\{t_\lambda : \lambda \in \Lambda\}$, there is a
homomorphism $\pi_t$ satisfying $\pi_t(s_\lambda) = t_\lambda$ for
all $\lambda \in \Lambda$.
By \cite[Proposition~3.5]{RSY1}, if $\mu,\nu \in \Lambda$, then
$s^*_\mu s_\nu = \sum_{\mu\alpha = \nu\beta, d(\mu\alpha) =
d(\mu) \vee d(\nu)} s_\alpha s^*_\beta$, and hence
(\cite[Remarks~3.8(1)]{RSY1}),
\begin{equation}\label{spanningset}
C^*(\Lambda) = \clsp\{s_\alpha s^*_\beta : s(\alpha) = s(\beta)\}.
\end{equation}

For the details of the next two paragraphs, see \cite[page~109]{RSY1}.

The universal property of $C^*(\Lambda)$ guarantees that there
is an action $\gamma : \T^k \to \Aut(C^*(\Lambda))$ satisfying
$\gamma_z(s_\lambda) = z^{d(\lambda)}
s_\lambda:=z_1^{d(\lambda)_1}\cdot z_2^{d(\lambda)_2}\cdots
z_k^{d(\lambda)_k}s_\lambda$ and hence $\gamma_z(p_v) = p_v$. We
denote the fixed point algebra for $\gamma$ by $F$, and $\Phi$
denotes the faithful conditional expectation $\Phi : C^*(\Lambda)
\to F$ determined by $\Phi(a) = \int_{\T^k} \gamma_z(a)\,d\mu(z)$.

We have $F = \clsp\{s_\mu s^*_\nu : d(\mu) = d(\nu), s(\mu) =
s(\nu)\}$ and $\Phi$ is determined by $\Phi(s_\mu s^*_\nu) =
\delta_{d(\mu),d(\nu)} s_\mu s^*_\nu$. For each $n \in
\N^k$, we write
$F_n := \clsp\{s_\mu s^*_\nu : d(\mu) = d(\nu), \mu,\nu \in
\Lambda^{\le n}, s(\mu) = s(\nu)\}$.
Then each $F_n$ is isomorphic to a direct sum of matrix algebras and
algebras of compact operators, and $F = \overline{\bigcup F_n}$ is
an AF algebra.

\section{$k$-graph traces and faithful traces on $C^*(\Lambda)$}\label{traces}
In this section we investigate conditions which give rise to
faithful traces on $C^*(\Lambda)$ for a locally convex locally
finite $k$-graph $\Lambda$. As with the $C^*$-algebras of directed graphs, necessary and
sufficient conditions for the existence of faithful traces on a
$k$-graph algebra are hard to come by.
 We denote by $A^+$ the positive cone in a
$C^*$-algebra $A$, and we use extended arithmetic on $[0,\infty]$ so that
 $0\times \infty=0$. From  \cite{PhR} we take the basic definition:

\begin{defn} A trace on a $C^*$-algebra $A$ is a map $\tau:A^+\to[0,\infty]$
satisfying

1) $\tau(a+b)=\tau(a)+\tau(b)$ for all $a,b\in A^+$

2) $\tau(\lambda a)=\lambda\tau(a)$ for all $a\in A^+$ and $\lambda\geq 0$

3) $\tau(a^*a)=\tau(aa^*)$ for all $a\in A$

We say: that $\tau$ is faithful if $\tau(a^*a)=0\Rightarrow a=0$; that $\tau$
 is semifinite if $\{a\in A^+:\tau(a)<\infty\}$ is norm dense in $A^+$ (or
that $\tau$ is densely defined); that $\tau$ is lower semicontinuous if
 whenever $a=\lim_{n\to\infty}a_n$ in norm in $A^+$ we have
$\tau(a)\leq\lim\inf_{n\to\infty}\tau(a_n)$.
\end{defn}
We may extend a (semifinite) trace $\tau$  by linearity to a
linear functional on  (a dense subspace of) $A$. Observe that the domain
of definition of a densely defined trace is a two-sided  ideal
$I_\tau\subset A$.
The proof of the following Lemma is identical to that of the analogous
result for graph algebras \cite[Lemma 3.2]{PRen}.
\begin{lemma}\label{finiteonfinite} Let $(\Lambda,d)$ be a row-finite
locally convex $k$-graph and let
$\tau:C^*(\Lambda)\to\C$ be a semifinite trace. Then the dense subalgebra
$$ A_c:={\rm span}\{s_\mu s_\nu^*:\mu,\nu\in \Lambda\}$$
is contained in the domain $I_\tau$ of $\tau$.
\end{lemma}

Recall from \cite{Si} that a loop with an entrance is a path
$\lambda\in\Lambda$ with $r(\lambda)=s(\lambda)$ such that
$d(\lambda)\geq e_i$ for some $1\leq i\leq k$, together with an
$e\in\Lambda^{e_i}$ with $r(e)=r(\lambda)$ but $\lambda(0,e_i)\neq
e$.

\begin{lemma}\label{necessary}  Let $(\Lambda,d)$ be a row-finite locally convex
$k$-graph.
\par\noindent {\bf (i)} If $C^*(\Lambda)$ has a faithful semifinite trace
then no loop
can have an entrance.
\par\noindent
{\bf (ii)} If $C^* (\Lambda)$ has a gauge-invariant, semifinite, lower
semicontinuous trace $\tau$ then
$\tau \circ \Phi = \tau$ and
$$
\tau(s_\mu s_\nu^*)=\delta_{\mu,\nu}\tau(p_{s(\mu)}).
$$

\noindent In particular, if $\tau|_{C^* (\{s_\mu s_\mu^* :
\mu \in \Lambda\})} = 0$ then $\tau = 0$.
\end{lemma}

\begin{figure}[ht]
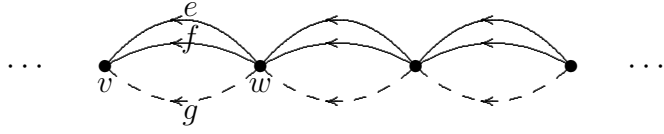

\[
\beginpicture
\setcoordinatesystem units <2.5em,2.5em>
\put{ } at 0 0.6
\put{ } at 0 -0.3
\put{$\cdots$} at -1 0
\put{$\bullet$} at 0 0
\put{$v$}[t] at 0 -0.15
\put{$\bullet$} at 2 0
\put{$w$}[t] at 2 -0.15
\put{$\bullet$} at 4 0
\put{$\bullet$} at 6 0
\put{$\cdots$} at 7 0
\setquadratic
\plot 0 0 1 0.3 2 0 /
\plot 0 0 1 0.6 2 0 /
\plot 2 0 3 0.3 4 0 /
\plot 2 0 3 0.6 4 0 /
\plot 4 0 5 0.3 6 0 /
\plot 4 0 5 0.6 6 0 /
\setdashes
\plot 0 0 1 -0.45 2 0 /
\plot 2 0 3 -0.45 4 0 /
\plot 4 0 5 -0.45 6 0 /
\setsolid
\arr (.86, 0.3)(.85,0.3)
\put{$f$}[B] at 1.1 0.25
\arr (.86, 0.6)(.85,0.6)
\put{$e$}[B] at 1.1 0.65
\arr (.86, -0.45)(.85,-0.45)
\put{$g$}[t] at 1.1 -0.5
\arr (2.86, 0.3)(2.85,0.3)
\arr (2.86, 0.6)(2.85,0.6)
\arr (2.86, -0.45)(2.85,-0.45)
\arr (4.86, 0.3)(4.85,0.3)
\arr (4.86, 0.6)(4.85,0.6)
\arr (4.86, -0.45)(4.85,-0.45)
\endpicture
\]
\caption{A $2$-graph whose $C^*$-algebra does not admit a faithful semifinite trace}
\label{fig:no trace}
\end{figure}

\begin{proof}
The entrance condition implies that $\lambda(0, e_i)$ and the entrance
$e$ are distinct paths of degree $e_i$ with the same range, and it follows
from (CK3)~and~(CK4) that
\[
s^*_\lambda s_\lambda = p_{s(\lambda)} = p_{r(\lambda)}
\ge s_\lambda s^*_\lambda + s_e s^*_e.
\]
If $\tau$ is a trace on $C^*(\Lambda)$, we therefore have
$\tau(s^*_\lambda s_\lambda) \ge \tau(s_\lambda s^*_\lambda)
+ \tau(s_e s^*_e)$, and it follows from Lemma \ref{finiteonfinite} and
the trace property
that $\tau(s^*_e s_e) = \tau(s_e s^*_e) = 0$. Theorem~3.15 of \cite{RSY1}
implies that $s^*_e s_e \not= 0$ so $\tau$ is not faithful.

The proof of the second part is the same as \cite[Lemma 3.3]{PRen},
but for clarity we remind the reader how the final statement arises.
If $\tau$ is gauge invariant we have
$$\tau(s_\mu s_\nu^*)=\tau(\gamma_z (s_\mu
s_\nu^*))=z^{d(\mu)-d(\nu)}\tau(s_\mu s_\nu^*)$$ for
all $z\in \T^k$. Hence $\tau(s_\mu s_\nu^*)$ is zero unless
$d(\mu) = d(\nu)$, and so $\tau=\tau\circ\Phi$.
Moreover if $d(\mu) = d(\nu)$, then using the trace property,
$$\tau(s_\mu
s_\nu^*)=\tau(s_\nu^* s_\mu)=\delta_{\nu,\mu}\tau(p_{s(\nu)})
=\delta_{\nu,\mu}\tau(s_\nu^* s_\nu).$$ This proves that if
$\tau|_{\lsp\{s_\mu s^*_\mu : \mu \in \Lambda\}} = 0$ then $\tau|{A_c} = 0$.
The details of extending this to the $C^*$-completion are as in \cite{PRen}.
\end{proof}
Whilst the condition that no loop has an entrance is necessary for the existence
of a faithful semifinite trace, it is not sufficient. For example, let $\Lambda$
be any $2$-graph whose skeleton is the one illustrated in Figure~\ref{fig:no trace}
(there are many allowable factorisation regimes to choose from).
Then $\Lambda$ is locally convex and locally finite, contains no sinks or sources,
and contains no cycles at all, so certainly no cycles with entrances, yet
$C^*(\Lambda)$ does not admit a faithful semifinite trace.
To see why
note that~(CK4) forces $s_g s^*_g = p_v = s_e s^*_e + s_f s^*_f$ so if $\tau$
is a trace on $C^*(\Lambda)$ then the trace property forces
\[
\tau(p_v) = \tau(s_g s^*_g) = \tau(p_w)
\quad\text{and}\quad
\tau(p_v) = \tau(s_e s^*_e) + \tau(s_f s^*_f) = 2\tau(p_w),
\]
and hence $\tau(p_w) = 0$.

The situation illustrated in Figure~\ref{fig:no trace} is more subtle than those which
can arise for graph $C^*$-algebras. However, as with directed graphs, the obstructions
to the existence of a faithful semifinite trace on a $k$-graph algebra can be expressed
most naturally for general $k$-graphs $\Lambda$ in terms of a
function $g_\tau : \Lambda^0 \to \R^+$ which arises naturally from each trace $\tau$
on $C^*(\Lambda)$.

\begin{lemma}\label{lem:trace gives graph trace}
Let $\Lambda$ be a locally convex row-finite $k$-graph, and
suppose that $\tau$ is a semifinite trace on $C^*(\Lambda)$. Then the
function $g_\tau : \Lambda^0 \to \R^+$ defined by $g_\tau(v) :=
\tau(p_v)$ satisfies $g_\tau(v) = \sum_{\lambda \in v\Lambda^{\le
n}} g_\tau(s(\lambda))$ for all $v \in \Lambda^0$ and $n \in
\N^k$.
\end{lemma}
\begin{proof}
Fix $v \in \Lambda^0$ and $n \in \N^k$. By~(CK4), we have $p_v =
\sum_{\lambda \in v\Lambda^{\le n}} s_\lambda s^*_\lambda$. Hence
\[
\tau(p_v) = \sum_{\lambda \in v\Lambda^{\le n}} \tau(s_\lambda
s^*_\lambda) = \sum_{\lambda \in v\Lambda^{\le n}}
\tau(s^*_\lambda s_\lambda) = \sum_{\lambda \in v\Lambda^{\le n}}
\tau(p_{s(\lambda)}),
\]
and the result follows from the definition of $g_\tau$.
\end{proof}

Returning to the example of Figure~\ref{fig:no trace} we can see that if
$\tau$ is a trace on $C^*(\Lambda)$ then $g_\tau(v)$ must simultaneously
be equal to $g_\tau(w)$ and $2g_\tau(w)$, forcing $g_\tau(w)$ and hence
$\tau(p_w)$ to be equal to zero.

Motivated by Lemma~\ref{lem:trace gives graph trace}, we make the
following definition (see \cite{H} for the origins of this definition):

\begin{defn}\label{dfn:graph trace}
Let $\Lambda$ be a locally convex row-finite $k$-graph. A function
$g : \Lambda^0 \to \R^+$ is called a \emph{$k$-graph trace} on
$\Lambda$ if it satisfies
\begin{equation}\label{eq:graph trace}
g(v) = \sum_{\lambda \in v\Lambda^{\le n}} g(s(\lambda))\text{ for
all $v \in \Lambda^0$ and $n \in \N^k$}.
\end{equation}
We say that $g$ is \emph{faithful} if $g(v) \not= 0$ for all $v
\in \Lambda^0$.
\end{defn}

\begin{rems}\label{rmk:const on ends}
Notice that if $x$ is an end of $\Lambda$, then
$x(0)\Lambda^{\le n} = \{x(0,n)\}$ for any $n \le d(x)$. It follows
that each $k$-graph trace on $\Lambda$ is constant on the vertices of $x$.
\end{rems}

We want to be able to construct semifinite lower semicontinuous gauge-invariant traces on
$C^*(\Lambda)$ from $k$-graph traces on $\Lambda$. The idea is to
use~\eqref{eq:graph trace} to define a trace on $C^*(\Lambda)$ by
$\tau_g\Big(\sum_{\mu, \nu \in F} a_{\mu,\nu} s_\mu s^*_\nu\Big)
= \sum_{\mu \in F} a_{\mu,\mu}g(s(\mu))$. There are two problems to overcome:
is $\tau_g$ well-defined in the first place, and when is $\tau_g$
faithful? To address these problems, we establish that there is a faithful
conditional expectation $\Psi$ on
$C^*(\Lambda)$ satisfying $\Psi(s_\mu s^*_\nu) =
\delta_{\mu,\nu} s_\mu s^*_\mu$. We would like to thank the referee
for pointing out the straightforward proof of this result appearing
below.

\begin{prop} \label{prp:diagonal map}
Let $\Lambda$ be a locally convex row-finite $k$-graph. There is a
faithful conditional expectation $\Psi : C^*(\Lambda) \to D :=
\clsp\{s_\lambda s^*_\lambda : \lambda \in \Lambda\}$ which
satisfies $\Psi(s_\lambda s^*_\mu) = \delta_{\lambda,\mu}
s_\lambda s^*_\lambda$ for all $\lambda,\mu \in \Lambda$.
\end{prop}
\begin{proof}
Averaging over the gauge action $\gamma$ gives a faithful conditional
expectation $\Phi^\gamma$ onto the fixed point algebra
$C^*(\Lambda)^\gamma = \clsp\{s_\mu s^*_\nu : d(\mu) = d(\nu)\}$.
For $q \in \NN^k$, $p \le q$ and $v \in \Lambda^0$, let
$P(q,p,v) := \{\lambda \in \Lambda^{\le q} : s(\lambda) = v, d(\lambda) = p\}$.
It is shown on page~109 of \cite{RSY1} that
$\mathcal{F}_{q,p}(v):= \clsp\{s_\mu s^*_\nu : \mu,\nu \in P(q,p,v)\}$,
is canonically isomorphic to $\mathcal{K}(\ell^2(P(q,p,v)))$, that for fixed
$q$, distinct $\mathcal{F}_{q,p}(v)$ are orthogonal, and that setting
$\mathcal{F}_q := \bigoplus_{p,v} \mathcal{F}_{q,p}(v)$, we have
\[
    C^*(\Lambda)^\gamma = \overline{\bigcup_{q \in \NN^k} \mathcal{F}_q}
\]
This shows that $C^*(\Lambda)^\gamma$ is an AF algebra with  maximal abelian
subalgebra $D$. Let $\Phi^D$ denote the canonical expectation
from $C^*(\Lambda)^\gamma$ onto its maximal abelian subalgebra. Then
$\Psi := \Phi^D \circ \Phi^\gamma$ is the desired expectation.
\end{proof}

\begin{prop}\label{trace=graphtrace} Let $\Lambda$ be a row-finite
locally convex $k$-graph.
Then there is a one-to-one correspondence between faithful graph
traces on $\Lambda$ and faithful, semifinite, lower semicontinuous,
gauge invariant traces on $C^*(\Lambda)$.
\end{prop}

\begin{proof}
Given a faithful $k$-graph trace, the existence of
$\Psi:C^*(\Lambda)\to\overline{span}\{S_\mu S_\mu^*\}$ given by
Proposition~\ref{prp:diagonal map} shows
that the functional $\tau_g:A_c\to\C$ defined by
$$\tau_g(S_\mu S_\nu^*):=\delta_{\mu,\nu}g(s(\mu))$$
is well-defined. As in \cite{PRen}, one checks that $\tau_g$ is a
gauge invariant
trace on $A_c$ and is faithful  because for
$a = \sum_{i=1}^n c_{\mu_i , \nu_i } S_{\mu_i} S_{\nu_i}^* \in
A_c$ we have  $a^* a \ge \sum_{i=1}^n \vert c_{\mu_i , \nu_i}
\vert^2 S_{\nu_i} S_{\nu_i}^*$, so
\begin{align}
\langle a , a \rangle_g &:= \tau_g ( a^* a ) = \tau_g \Big(\sum_{i=1}^n
\vert c_{\mu_i , \nu_i} \vert^2 S_{\nu_i} S_{\nu_i}^*\Big) =
\sum_{i=1}^n \vert c_{\mu_i , \nu_i} \vert^2 \tau_g ( S_{\nu_i}
S_{\nu_i}^* ) = \sum_{i=1}^n \vert c_{\mu_i , \nu_i } \vert^2 g (
s ( \nu_i ) ) > 0 .
\end{align}
by definition of $\tau_g$.

Then $\la a,b\ra_g=\tau_g(b^*a)$ defines a positive definite inner product
on $A_c$ which makes it a Hilbert algebra (that the left regular
representation of $A_c$ is nondegenerate follows from $A_c^2=A_c$).
The rest of the proof is the same as \cite[Proposition 3.9]{PRen},
except that we use the gauge invariant uniqueness theorem for
$k$-graphs, \cite[Theorem~4.1]{RSY1}, to show that we obtain
faithful representation of
$A=C^*(\Lambda)$ on the Hilbert space completion of $A_c$.
\end{proof}

The results of sections \ref{spectrip}--\ref{index} below apply to
any $k$-graph which admits a faithful graph trace. It is therefore
important to establish that there is a substantial class of
$k$-graphs for any $k$ with this property, and how large the class
is. This is in general a difficult question. The results of
\cite{PRRS} show there is a substantial class of $2$-graphs
admitting faithful graph traces, but for arbitrary $k$ there are, to
our knowledge, no definitive results. In Appendix~\ref{appendix} we
establish some necessary conditions and one sufficient condition,
which will prove useful in section~\ref{index}.

\section{Semifinite Spectral Triples}\label{spectrip}

We begin with some semifinite versions of standard definitions and
results. Let $\tau$ be a fixed faithful, normal, semifinite trace
on the von Neumann algebra ${\mathcal N}$. Let ${\mathcal
K}_{\mathcal N }$ be the $\tau$-compact operators in ${\mathcal
N}$ (that is the norm closed ideal generated by the projections
$E\in\mathcal N$ with $\tau(E)<\infty$).

\begin{defn} A semifinite
spectral triple $(\A,\HH,\D)$ relative to $(\cn,\tau)$ consists
of a Hilbert space $\HH$, a
$*$-algebra $\A\subset \cn$ where $\cn$ is a semifinite von
Neumann algebra acting on $\HH$, and a densely defined unbounded
self-adjoint operator $\D$ affiliated to $\cn$ such that

1) $[\D,a]$ is densely defined and extends to a bounded operator
in $\cn$ for all $a\in\A$

2) $a(\lambda-\D)^{-1}\in\K_\cn$ for all $\lambda\not\in{\R}$ and
all $a\in\A$

3) The triple is said to be even if there is $\Gamma\in\cn$ such
that $\Gamma^*=\Gamma$, $\Gamma^2=1$, $\D\Gamma+\Gamma\D=0$,
and $a\Gamma=\Gamma a$ for all $a\in\A$. Otherwise it is odd.
\end{defn}

\begin{defn}\label{qck} A semifinite spectral triple $(\A,\HH,\D)$ is $QC^k$ for $k\geq 1$
($Q$ for quantum) if for all $a\in\A$ the operators $a$ and
$[\D,a]$ are in the domain of $\delta^k$, where
$\delta(T)=[\dd,T]$ is the partial derivation on $\cn$ defined by
$\dd$. We say that $(\A,\HH,\D)$ is $QC^\infty$ if it is $QC^k$
for all $k\geq 1$.
\end{defn}

{\bf Note}. The notation is meant to be analogous to the classical
case, but we introduce the $Q$ so that there is no confusion
between quantum differentiability of $a\in\A$ and classical
differentiability of functions.

\noindent{\bf Remarks concerning derivations and commutators}.  By
partial derivation we mean that $\delta$ is defined on some
subalgebra of $\cn$ which need not be (weakly) dense in $\cn$.
More precisely, $\mbox{dom}\ \delta=\{T\in\cn:\delta(T)\mbox{ is
bounded}\}$. We also note that if $T\in{\mathcal N}$, one can show
that $[\dd,T]$ is bounded if and only if $[(1+\D^2)^{1/2},T]$ is
bounded, by using the functional calculus to show that
$\dd-(1+\D^2)^{1/2}$ extends to a bounded operator in $\cn$. In
fact, writing $\dd_1=(1+\D^2)^{1/2}$ and $\delta_1(T)=[\dd_1,T]$
we have \ben \mbox{dom}\ \delta^n=\mbox{dom}\ \delta_1^n\ \ \ \
\mbox{for all}\ n.\een We also observe that if $T\in\cn$ and
$[\D,T]$ is bounded, then $[\D,T]\in\cn$. Similar comments apply
to $[\dd,T]$, $[(1+\D^2)^{1/2},T]$. The proofs of these statements
can be found in \cite{CPRS2}.
\begin{defn}A $*$-algebra $\A$ is smooth if it is Fr\'{e}chet
and $*$-isomorphic to a proper dense subalgebra $i(\A)$ of a
$C^*$-algebra $A$ which is stable under the holomorphic functional
calculus.\end{defn} Thus saying that $\A$ is \emph{smooth} means
that $\A$ is Fr\'{e}chet and a pre-$C^*$-algebra. Asking for
$i(\A)$ to be a {\it proper} dense subalgebra of $A$ immediately
implies that the Fr\'{e}chet topology of $\A$ is finer than the
$C^*$-topology of $A$ (since Fr\'{e}chet means locally convex,
metrizable and complete.) We will sometimes speak of
$\overline{\A}=A$, particularly when $\A$ is represented on
Hilbert space and the norm closure $\overline{\A}$ is unambiguous.
At other times we regard $i:\A\hookrightarrow A$ as an embedding
of $\A$ in a $C^*$-algebra. We will use both points of view.

It has been shown that if $\A$ is smooth in $A$ then $M_n(\A)$ is
smooth in $M_n(A)$, \cite{GVF,LBS}. This ensures that the
$K$-theories of the two algebras are isomorphic, the isomorphism
being induced by the inclusion map $i$. This definition ensures
that a smooth algebra is a `good' algebra, \cite{GVF}, so these
algebras have a sensible spectral theory which agrees with that
defined using the $C^*$-closure, and the group of invertibles is
open.

Stability under the holomorphic functional calculus extends to
nonunital algebras, since the spectrum of an element in a
nonunital algebra is defined to be the spectrum of this element in
the  `one-point' unitization, though we must of course restrict to
functions satisfying $f(0)=0$. Likewise, the definition of a
Fr\'{e}chet algebra does not require a unit. The point of contact
between smooth algebras and $QC^\infty$ spectral triples is the
following Lemma, proved in \cite{R1}.

\begin{lemma}\label{smo} If $(\A,\HH,\D)$ is a $QC^\infty$ spectral triple, then
$(\A_\delta,\HH,\D)$ is also a $QC^\infty$ spectral triple, where
$\A_\delta$ is the completion of $\A$ in the locally convex
topology determined by the seminorms \ben
q_{n,i}(a)=\|\delta^nd^i(a)\|,\ \ n\geq 0,\ i=0,1,\een where
$d(a)=[\D,a]$. Moreover, $\A_\delta$ is a smooth algebra.
\end{lemma}

We call the topology on $\A$ determined by the seminorms $q_{ni}$
of Lemma \ref{smo} the $\delta$-topology.

Whilst smoothness does not depend on whether $\A$ is unital or
not, many analytical problems arise because of the lack of a unit.
As in \cite{R1,R2,GGISV}, we make two definitions to address these
issues.

\begin{defn} An algebra $\A$ has local units if for every finite subset of
elements $\{a_i\}_{i=1}^n\subset\A$, there exists $\phi\in\A$ such
that for each $i$ \ben \phi a_i= a_i\phi=a_i.\een
\end{defn}

\begin{defn}
Let $\A$ be a Fr\'{e}chet algebra and $\A_c\subseteq\A$ be a dense
subalgebra with local units. Then we call  $\A$ a quasi-local
algebra (when $\A_c$ is understood.) If $\A_c$ is a dense ideal
with local units, we call $\A_c\subset\A$ local.
\end{defn}

Quasi-local algebras have an approximate unit $\{\phi_n\}_{n\geq
1}\subset\A_c$ such that $\phi_{n+1}\phi_n=\phi_n$, \cite{R1}.

{\bf Example} For a $k$-graph $C^*$-algebra $A=C^*(\Lambda)$,
Equation (\ref{spanningset}) shows that
$$ A_c=\mbox{span}\{S_\mu S_\nu^*:\mu,\nu\in E^*\ \mbox{and}\
s(\mu)=s(\nu)\}$$ is a dense subalgebra. It has local units
because
$$ p_{v}S_\mu S_\nu^*=\left\{\begin{array}{ll} S_\mu S_\nu^* &
v=r(\mu)\\ 0 & \mbox{otherwise}\end{array}\right..$$ Similar
comments apply to right multiplication by $p_{r(\nu)}$. By summing
the source and range projections (without repetitions) of all
$S_{\mu_i}S_{\nu_i}^*$ appearing in a finite sum
$$ a=\sum_ic_{\mu_i,\nu_i}S_{\mu_i}S_{\nu_i}^*$$
we obtain a local unit for $a\in A_c$. By repeating this process
for any finite collection of such $a\in A_c$ we see that $A_c$ has
local units.

We also require that when we have a spectral triple the operator
$\D$ is compatible with the quasi-local structure of the algebra,
in the following sense.

\begin{defn} If $(\A,\HH,\D)$ is a spectral triple, then we define $\Omega^*_\D(\A)$
to be the algebra generated by $\A$ and $[\D,\A]$.
\end{defn}

\begin{defn}\label{lst} A local spectral triple $(\A,\HH,\D)$ is a
spectral triple with $\A$ quasi-local such that there exists an
approximate unit $\{\phi_n\}\subset\A_c$ for $\A$ satisfying \ben
\Omega^*_\D(\A_c)=\bigcup_n\Omega^*_\D(\A)_n,\een \ben
\Omega^*_\D(\A)_n=\{\omega\in\Omega^*_\D(\A):\phi_n\omega=\omega\phi_n=\omega\}.\een
\end{defn}

{\bf Remark} A local spectral triple has a local approximate unit
$\{\phi_n\}_{n\geq 1}\subset\A_c$  such that, \cite{R2},
$\phi_{n+1}\phi_n=\phi_n\phi_{n+1}=\phi_n$ and
$\phi_{n+1}[\D,\phi_n]=[\D,\phi_n]\phi_{n+1}=[\D,\phi_n]$. This is
the crucial property we require to prove our summability
results for nonunital spectral triples, to which we now turn.

\subsection{Summability}
In the following, let $\mathcal N$ be a semifinite von Neumann
algebra with faithful normal trace $\tau$. Recall from \cite{FK}
that if $S\in\mathcal N$, the \emph{$t^{\rm th}$ generalized singular
value} of S for each real $t>0$ is given by
$$\mu_t(S)=\inf\{\|SE\|\ : \ E \mbox{ is a projection in }
{\mathcal N} \mbox { with } \tau(1-E)\leq t\}.$$

The ideal $\LL^1({\mathcal N})$ consists of those operators $T\in
{\mathcal N}$ such that $\|T\|_1:=\tau( |T|)<\infty$ where
$|T|=\sqrt{T^*T}$. In the Type~I setting this is the usual trace
class ideal. We will simply write $\LL^1$ for this ideal in order
to simplify the notation, and denote the norm on $\LL^1$ by
$\|\cdot\|_1$. An alternative definition in terms of singular
values is that $T\in\LL^1$ if $\|T\|_1:=\int_0^\infty \mu_t(T) dt
<\infty.$

Note that in the case where ${\mathcal N}\neq{\mathcal
B}({\mathcal H})$, $\LL^1$ is not complete in this norm but it is
complete in the norm $\|\cdot\|_1 + \|\cdot\|_\infty$. (where
$\|\cdot\|_\infty$ is the uniform norm). Another important ideal for
us is the domain of the Dixmier trace:
$${\mathcal L}^{(1,\infty)}({\mathcal N})=
\left\{T\in{\mathcal N}\ : \Vert T\Vert_{_{{\mathcal
L}^{(1,\infty)}}} :=   \sup_{t> 0}
\frac{1}{\log(1+t)}\int_0^t\mu_s(T)ds<\infty\right\}.$$

There are related ideals for $p>1$: to describe them first set
\[
\psi_p(t)=\begin{cases} t & \mbox{for } 0\leq t\leq 1 \\
                        t^{1-\frac{1}{p}} & \mbox{for } 1\leq t.
\end{cases}
\]
Then define
$${\mathcal L}^{(p,\infty)}({\mathcal N})
=\left\{T\in{\mathcal N}\ : \Vert T\Vert_{_{{\mathcal
L}^{(p,\infty)}}} := \sup_{t> 0}
\frac{1}{\psi_p(t)}\int_0^t\mu_s(T)ds<\infty\right\}.$$ For $p>1$
there is also the equivalent definition
$${\mathcal L}^{(p,\infty)}({\mathcal N})
=\left\{T\in{\mathcal N}\ : \sup_{t> 0}
\frac{t}{\psi_p(t)}\mu_t(T)<\infty\right\}.$$ If
$T\in\LL^{(p,\infty)}(\cn)$, then $T^p\in\LL^{(1,\infty)}(\cn)$.

We will suppress the $({\mathcal N})$ in our notation for these
ideals, as $\cn$ will always be clear from context. The reader
should note that ${\mathcal L}^{(1,\infty)}$ is often taken to
mean an ideal in the algebra $\widetilde{\mathcal N}$ of
$\tau$-measurable operators affiliated to ${\mathcal N}$. Our
notation is however consistent with that of \cite{C} in the
special case ${\mathcal N}={\mathcal B}({\mathcal H})$. With this
convention the ideal of $\tau$-compact operators, ${\mathcal
  K}({\mathcal N})$,
consists of those $T\in{\mathcal N}$ (as opposed to
$\widetilde{\mathcal N}$) such that \ben \mu_\infty(T):=\lim
_{t\to \infty}\mu_t(T)  = 0.\een

\begin{defn}\label{summable} A semifinite local spectral triple is
$(k,\infty)$-summable if \ben
a(\D-\lambda)^{-1}\in\LL^{(k,\infty)}\ \ \ \mbox{for all}\
a\in\A_c,\ \ \lambda\in\C\setminus\R.\een
\end{defn}

{\bf Remark} If $\A$ is unital, $\ker\D$ is $\tau$-finite
dimensional. Note that the summability requirements are only for
$a\in\A_c$. We do not assume that elements of the algebra $\A$ are
all integrable in the nonunital case. Strictly speaking, this
definition describes {\em local} $(k,\infty)$-summability, however
we use the terminology $(k,\infty)$-summable to be consistent with
the unital case.

We need to briefly discuss the Dixmier trace, but fortunately we
will usually be applying it in reasonably simple situations. For
more information on semifinite Dixmier traces, see \cite{CPS2}.
For $T\in\LL^{(1,\infty)}$, $T\geq 0$, the function
\ben
F_T : t \mapsto \frac{1}{\log(1+t)}\int_0^t\mu_s(T)ds
\een
is bounded. For
certain generalised limits $\omega\in L^\infty(\R_*^+)^*$, we
obtain a positive functional on $\LL^{(1,\infty)}$ by setting
$$ \tau_\omega(T)=\omega(F_T).$$
 This is the
Dixmier trace associated to the semifinite normal trace $\tau$,
denoted $\tau_\omega$, and we extend it to all of
$\LL^{(1,\infty)}$ by linearity, where of course it is a trace.
The Dixmier trace $\tau_\omega$ is defined on the ideal
$\LL^{(1,\infty)}$, and vanishes on the ideal of trace class
operators. Whenever the function $F_T$ has a limit at infinity,
all Dixmier traces return the value of the limit. We denote the
common value of all Dixmier traces on measurable operators by
$\bigintcross$. So if $T\in\LL^{(1,\infty)}$ is measurable, for
any allowed functional $\omega\in L^\infty(\R_*^+)^*$ we have
$$\tau_\omega(T)=\omega(F_T)=\bigintcross T.$$

{\bf Example} The Dirac operator on the $k$-torus. Let $\gamma^j$,
$j=1,\dots,k$, be generators of the Clifford algebra of $\R^k$
with the usual Euclidean inner product. Form the Dirac operator on
spinors $\D=\sum_{j=1}^k\gamma^j\frac{\p}{\p\theta^j}$, which acts
on $L^2(\T^k)\otimes\C^{2^{[k/2]}}$, and for $n \in \Z^k$, let
$n^2 \in \N$ denote the sum $n^2 = \sum^k_{i=1} n_i^2$ of the squares
of the coordinates of $n$. Then it is well known that
the spectrum of $\D^2$ consists of eigenvalues $\{n^2\in\N\}$,
where each $n\in\Z^k$ is counted once.  A careful calculation
taking account of the multiplicities, \cite{La}, shows that using
the standard operator trace, the function $F_{(1+\D^2)^{-k/2}}$ is
$$ \frac{1}{\log (|\{n:|n|\leq N\}|)}\sum_{|n|=0}^N(1+n^2)^{-k/2}=
\frac{2^{[k/2]}vol(S^{k-1})}{k\log
N}\sum_{m=0}^N(1+m^2)^{-1/2}+o(1)$$ and this is bounded. Hence
$(1+\D^2)^{-k/2}\in\LL^{(1,\infty)}$ and
$$\mbox{Trace}_\omega((1+\D^2)^{-k/2})=\bigintcross(1+\D^2)^{-k/2}
=\frac{2^{[k/2]}vol(S^{k-1})}{k}=\frac{2^{[k/2]}vol(S^{k-1})}
{(2\pi)^kk}vol(\T^k).$$

Numerous properties of local algebras are established in \cite{R1,R2}.
The introduction of quasi-local algebras in \cite{GGISV} led to a
review of the validity of many of these results for quasi-local
algebras. Most of the summability results of \cite{R2} are valid
in the quasi-local setting.  In addition, the summability results
of \cite{R2} are also valid for general semifinite spectral
triples since they rely only on properties of the ideals
$\LL^{(p,\infty)}$, $p\geq 1$, \cite{C,CPS2}, and the trace
property. We quote the version of the summability results from
\cite{R2} that we require below.

\begin{prop}[\cite{R2}]\label{wellbehaved} Let $(\A,\HH,\D)$ be a $QC^\infty$, local
$(k,\infty)$-summable semifinite spectral triple. Let $T\in\cn$
satisfy $T\phi=\phi T=T$ for some $\phi\in\A_c$. Then \ben
T(1+\D^2)^{-k/2}\in\LL^{(1,\infty)}.\een For $Re(s)>k$,
$T(1+\D^2)^{-s/2}$ is trace class. If the limit \be \lim_{s\to
k/2^+}(s-k/2)\tau(T(1+\D^2)^{-s})\label{mumbo}\ee exists, then it
is equal to \ben \frac{k}{2}\bigintcross T(1+\D^2)^{-k/2}.\een In
addition, for any Dixmier trace $\tau_\omega$, the function \ben
a\mapsto \tau_\omega(a(1+\D^2)^{-k/2})\een defines a trace on
$\A_c\subset\A$.
\end{prop}

\section{Constructing a $C^*$-module and a Kasparov module}\label{firstconstruction}
Let $A=C^*(\Lambda)$ where $\Lambda$ is a locally finite
locally convex $k$-graph. Let $F=C^*(\Lambda)^\gamma$ be the
fixed point subalgebra for the gauge action. Finally, let
$A_c = \lsp\{s_\mu s^*_\nu : \mu,\nu \in \Lambda\}$ and let
$F_c = \lsp\{s_\mu s^*_\nu : d(\mu) = d(\nu)\} = F \cap A_c$ so
that $A$ and $F$ are the $C^*$-completions of $A_c$ and $F_c$. Note
that the expectation $\Phi : A \to F$ outlined at the end of
Section~\ref{kgraph} restricts to an expectation, also denoted $\Phi$
of $A_c$ onto $F_c$.

For $q \in \Q$, we denote by $[q]$ the integer part $\max\{n \in \ZZ :
n \le q\}$ of $q$.
We make $A^{2^{[k/2]}}=\C^{2^{[k/2]}}\otimes A$ a right
inner product-$F$-module. The right action of $F$ on $A$ is by right
multiplication.
The inner product is defined by $$
(x|y)_R:=\sum_{j=1}^{2^{[k/2]}}\Phi(x_j^*y_j)\in F.$$ It is simple to
check the
requirements that $(\cdot|\cdot)_R$ defines an $F$-valued inner
product on $A^{2^{[k/2]}}$. The requirement $(x|x)_R=0\Rightarrow x=0$
follows
from the faithfulness of $\Phi$.

\begin{defn}\label{Fmod} Define $X$ to be the completion of $A^{2^{[k/2]}}$
to a $C^*$-module over $F$ for the $C^*$-module norm
$$\Vert x\Vert_X^2:=\Vert(x|x)_R\Vert_A=\Vert(x|x)_R\Vert_F=
\Vert \sum_{i=1}^{2^{[k/2]}}\Phi(x_i^*x_i)\Vert_F.$$
Define $X_c$
to be the pre-$C^*$-module over $F_c$ with linear space
$A_c^{2^{[k/2]}}$ and the inner product $(\cdot|\cdot)_R$.
\end{defn}

{\bf Remark} Typically, the action of $F$ does not map $X_c$ to
itself, so we may only consider $X_c$ as an $F_c$ module. This is
a reflection of the fact that $F_c$ and $A_c$ are quasilocal not
local.

{\bf Remark} Frequently we  will define an operator $T$ on the $F$
module $A$, and implicitly extend $T$  to $X$ by $\id_{2^{[k/2]}} \otimes T$,
where $\id_{2^{[k/2]}}$ is the identity operator
in the matrix algebra $M_{2^{[k/2]}}(\C)$.

{\bf Remark} There is an irreducible representation $\gamma$ of
the complex Clifford algebra $\C\textit{liff}_k=\textit{Cliff}(\C^k)$
on $\C^{2^{[k/2]}}$, and tensoring this representation by the identity
map on $A$,
this extends to a representation on $X$ as
adjointable operators. We employ the convention that
$$
\gamma^l\gamma^j+\gamma^j\gamma^l
:=\gamma(e^l)\gamma(e^j)+\gamma(e^j)\gamma(e^l)=-2\delta^{lj}\id_{2^{[k/2]}}.
$$
When $k$ is even the operator
$\omega_\C:=i^{[(k+1)/2]}\gamma^1\cdots\gamma^k$ is self-adjoint,
has $\omega^2_\C=\id_{2^{[k/2]}}$ and
$\gamma^j\omega_\C=-\omega_\C\gamma^j$ for $j=1,\dots,k$. When
$k$ is odd, $\omega_\C$ is central in the Clifford algebra, and
we choose the representation with $\omega_\C=1$.

The map $a \mapsto 1_{2^{[k/2]}} \otimes a$ is an
isometric inclusion of $A$ into $\C^{2^{[k/2]}}\otimes A
= A^{2^{[k/2]}}$, which in turn is dense in $X$
by definition. The inclusion $\iota : A \to X$ is continuous since
$$\Vert a\Vert_X^2=\Vert\Phi(a^*a)\Vert_F\leq\Vert
a^*a\Vert_A=\Vert a\Vert^2_A.$$ We can also define the gauge
action $\gamma$ on $A\subset X$, and as
\bean\Vert\gamma_z(a)\Vert^2_X&=&\Vert\Phi((\gamma_z(a))^*(\gamma_z(a)))\Vert_F
=\Vert\Phi(\gamma_z(a^*)\gamma_z(a))\Vert_F\nno&=&
\Vert\Phi(\gamma_z(a^*a))\Vert_F =\Vert\Phi(a^*a)\Vert_F=\Vert
a\Vert^2_X,\eean for each $z\in \T^k$, the action of $\gamma_z$ is
isometric on $A\subset X$ and so extends to a unitary $U_z$ on
$X$. This unitary is $F$-linear and adjointable, and we obtain a
strongly continuous action of $\T^k$ on $X$, which we still denote
by $\gamma$.

For each $n\in\Z^k$, define an operator $\Phi_n$ on $X$  by
$$\Phi_n(x)=\frac{1}{(2\pi)^k}\int_{\T^k}z^{-n}\gamma_z(x)d^k\theta,\ \
z_j=e^{i\theta_j},\ \ x\in X.$$
Observe that on generators we have
\begin{equation}\Phi_n(S_\al
S_\beta^*)=\left\{\begin{array}{lr}S_\al S_\beta^* & \ \
d(\al)-d(\beta)=n\\0 & \ \ d(\al)-d(\beta)\neq
n\end{array}\right..\label{nthproj}\end{equation}

{\bf Remark} If $(\Lambda,d)$ is a finite $k$-graph with no cycles,
then for $n$ sufficiently large there are no paths of
degree $n$ and so $\Phi_n=0$. This will obviously simplify many of
the convergence issues below.

The proof of the following Lemma is identical to that of
\cite[Lemma~4.2]{PRen}.

\begin{lemma}\label{phiendo}
The operators $\Phi_n$ are adjointable endomorphisms of the $F$-module $X$ such that $\Phi_n^*=\Phi_n=\Phi_n^2$ and
$\Phi_n\Phi_m=\delta_{n,m}\Phi_n$. For each subset $K\subset\Z^k$, the sum $\sum_{n\in K}\Phi_n$ converges strictly to a projection $\Phi_K$ in the
endomorphism algebra. Moreover, the projection $\Phi_{\ZZ^k}$ corresponding to $K = \ZZ^k$ is the identity operator on $X$.
\end{lemma}

\begin{cor}\label{gradedsum} Let $x\in X$. Then with $x_n=\Phi_nx$ the sum
$\sum_{n\in \Z^k}x_n$ converges in $X$ to $x$.
\end{cor}

\subsection{The Kasparov Module}\label{CstarDeeee}

As we did in Section~\ref{spectrip}, for
$n\in\Z^k$, we write $n^2=\sum_{j=1}^kn_j^2$ and $|n|=\sqrt{n^2}$.

 The theory of unbounded operators on
$C^*$-modules that we require is all contained in Lance's book,
\cite[Chapters 9,10]{L}. We quote the following definitions
(adapted to our situation).

\begin{defn} Let $Y$ be a right $C^*$-$B$-module. A densely defined
 unbounded operator $\D:{\rm dom}\ \D\subset Y\to Y$ is a
 $B$-linear operator defined on a dense $B$-submodule
 ${\rm dom}\ \D\subset Y$. The operator $\D$ is closed
 if the graph
 $$ G(\D)=\{(x,\D x)_R:x\in{\rm dom}\ \D\}$$
 is a closed submodule of $Y\oplus Y$.
\end{defn}

 Given a densely defined
 unbounded operator $\D:\mbox{dom}\ \D\subset Y\to Y$, define
 a submodule
 $$\mbox{dom}\ \D^*:=\{y\in Y:\exists z\in Y\ \mbox{such that}\
 \forall x\in\mbox{dom}\ \D, (\D x|y)_R=(x|z)_R\}.$$
 Then for $y\in \mbox{dom}\ \D^*$ define $\D^*y=z$. Given $y\in\mbox{dom}\ \D^*$,
 the element $z$ is unique, so $\D^*$ is well-defined, and
 moreover is closed.

 \begin{defn} Let $Y$ be a right $C^*$-$B$-module. A densely defined unbounded
 operator $\D:{\rm dom}\ \D\subset Y\to Y$ is symmetric if for all
 $x,y\in{\rm dom}\ \D$
 $$ (\D x|y)_R=(x|\D y)_R.$$
 A symmetric operator $\D$ is self-adjoint if
 ${\rm dom}\ \D={\rm dom}\ \D^*$ (and so $\D$ is necessarily
 closed). A densely defined unbounded operator $\D$ is regular if
 $\D$ is closed, $\D^*$ is densely defined, and $(1+\D^*\D)$ has
 dense range.
 \end{defn}

 The extra requirement of regularity is necessary in the
 $C^*$-module context for the continuous functional calculus,
 and is not automatically satisfied,
 \cite[Chapter 9]{L}.

With these definitions in hand, we return to our $C^*$-module $X$. The
proof of the following Proposition is an exact analogue of
\cite[Proposition 4.6]{PRen}.

\begin{prop}\label{CstarDee} Let $X$ be the right $C^*$-$F$-module of
Definition
 \ref{Fmod}.  Define $X_\D\subset X$ to be the linear space
$$X_\D=\{x=\sum_{n\in\Z^k}x_n:\Vert\sum_{n\in\Z^k}n^2(x_n|x_n)_R\Vert<\infty\}.$$
For $x\in X_\D$ define
$$\D
x=\sum_{n\in\Z^k}\gamma(in)x_n=i\sum_{n\in\Z^k}\sum_{j=1}^k\gamma^jn_jx_n.$$
Then $\D:X_\D\to X$
is self-adjoint and regular.
\end{prop}

{\bf Remark} For $n \in \ZZ^k$, the restriction of the map $\D$ to $\Phi_n X$ implements Clifford multiplication by the vector $in \in \C^k$. Any
$S_\al S_\beta^*\in A_c$ is in $X_\D$ and $$ \D S_\al S_\beta^*=i\sum_{j=1}^k\gamma^j(d(\al)_j-d(\beta)_j)S_\al S_\beta^*$$ as the reader will easily
verify. Thus we have $$\D^2\Phi_nx=\sum_{j=1}^kn_j^2\Phi_nx=n^2\Phi_nx.$$

There is a continuous functional calculus for self-adjoint regular
operators, \cite[Theorem 10.9]{L}, and we use this to obtain
spectral projections for $\D^2$ at the $C^*$-module level. Let
$f_m\in C_c({\R})$ be $1$ in a small neighbourhood of $m\in{\Z}$
and zero on $(-\infty,m-1/2]\cup[m+1/2,\infty)$. Then it is clear
that
$$ \sum_{n\in\Z^k,\ n^2=m}\Phi_n=f_m(\D^2).$$

The next Lemma is the first place where we need our $k$-graph to be
locally finite and have no sinks. It is also the point where the
generalisation from the graph case differs the most.

\begin{lemma}\label{finrank} Assume that the $k$-graph $(\Lambda,d)$ is
locally finite, locally convex and has no sinks. For all $a\in
A$ and $n\in\Z^k$, $a\Phi_n\in End^0_F(X)$, the
 compact  endomorphisms of the right $F$-module $X$. If $a\in A_c$ then
$a\Phi_n$ is finite rank.
\end{lemma}

\begin{rems} If we were employing the $A$-valued inner product on $X$, then
each $a\in A$ would be compact, and Lemma~\ref{finrank} would be an
immediate corollary of the fact that the compacts form an ideal.
However, with our choice of inner product, with values in $F$, {\em
no} $a\in A$ acts as a compact endomorphism, except in some extreme
examples.
\end{rems}

\begin{proof} Let $n\in\Z^k$ and write $n=n_1+n_2$ with $n_1\geq 0$
and $n_2<0$. We will see that the precise choice of $n_1,n_2$ is
largely irrelevant. For $v\in\Lambda^0$, let $|v|_n$ denote the number
of paths $\rho\in\Lambda$ with $d(\rho)=n$ and $s(\rho)=v$,
i.e. $|v|_n=|\Lambda^nv|$. Since $\Lambda$ has no sinks and is locally
finite, for all $n$ and all $v$ we have $0<|v|_n<\infty$ .

Now define, for $n=n_1+n_2$ as above,
$$T_{v,n_1,n_2}=
\sum_{d(\al)=n_1,d(\beta)=-n_2,s(\al)=s(\beta),r(\al)=v}
\frac{1}{|s(\beta)|_{-n_2}}
\Theta^R_{S_\al S_\beta^*,S_\al S_\beta^*},$$
where for $x,y,z\in X$
$$\Theta^R_{x,y}z:=x(y|z)_R,$$
defines a rank one operator. Observe that since $\Lambda$ is locally finite
this is a finite sum of rank one operators and so finite rank.
We claim that $T_{v,n_1,n_2}=p_v\Phi_n$.   It suffices to prove that
the difference $p_v\Phi_n-T_{v,n_1,n_2}$ vanishes on $X_c\subset X$. That is,
we just need to show that $(p_v \Phi_n - T_{v, n_1, n_2})S_\mu S^*_\nu = 0$
for all $\mu,\nu$. So first we compute, with $q=d(\al)\vee d(\mu)$,
$$S_\al^*S_\mu=\sum_{\al\s=\mu\rho,\al\s\in\Lambda^q}S_\s
S_\rho^*$$
by \cite[Proposition~3.5 and Remarks~3.8(2)]{RSY1}. Next consider
$$\Phi(S_\beta S_\al^*S_\mu
S_\nu^*)=\Phi(S_\beta\sum_{\al\s=\mu\rho,\al\s\in\Lambda^q}S_\s
S_\rho^*S_\nu^*).$$
This is zero unless $d(\beta)+d(\s)-d(\rho)-d(\nu)=0$. Now
$d(\s)-d(\rho)=d(\mu)-d(\al)$ so
$$\Phi(S_\beta S_\al^*S_\mu
S_\nu^*)=\delta_{d(\mu)-d(\nu),d(\al)-d(\beta)}
S_\beta\sum_{\al\s=\mu\rho,\al\s\in\Lambda^q}S_\s
S_\rho^*S_\nu^*.$$
Of course, $d(\al)-d(\beta)=n$. Since each
$S_\beta^*S_\beta=p_{s(\beta)}$, we can perform the sum over $\beta$:
\bean\sum_{\al,\beta}\frac{1}{|s(\beta)|_{-n_2}}p_v\Theta_{S_\al
S_\beta^*,S_\al S_\beta^*}S_\mu S_\nu^*&=&
\sum_{\al,\beta}
\frac{1}{|s(\beta)|_{-n_2}}\delta_{d(\mu)-d(\nu),n}
p_vS_\al S_\beta^*S_\beta\sum_{\al\s=\mu\rho,\al\s\in\Lambda^q}S_\s
S_\rho^*S_\nu^*\nno
&=&\sum_\al\delta_{d(\mu)-d(\nu),n}
p_vS_\al \sum_{\al\s=\mu\rho,\al\s\in\Lambda^q}S_\s
S_\rho^*S_\nu^*\nno
&=&\sum_\al\delta_{d(\mu)-d(\nu),n}
p_v \sum_{\al\s=\mu\rho,\al\s\in\Lambda^q}S_\mu S_\rho
S_\rho^*S_\nu^*.\eean
If we suppose that a given $\al$ has no common extensions with $\mu$,
then this particular term in the sum contributes zero. Summing over
all $\al$ (of fixed length $n_1$) with common extensions with $\mu$
yields
$$\sum_\al\delta_{d(\mu)-d(\nu),n}
p_v \sum_{\al\s=\mu\rho,\al\s\in\Lambda^q}S_\mu S_\rho
S_\rho^*S_\nu^*=p_v\sum_{\rho\in\Lambda^{q - d(\mu)}, r(\rho)=s(\mu)}S_\mu
S_\rho S_\rho^* S_\nu^*=S_\mu S_\nu^*.$$
Hence we conclude that
\bean \sum_{\al,\beta}\frac{1}{|s(\beta)|_{-n_2}}p_v\Theta_{S_\al
S_\beta^*,S_\al S_\beta^*}S_\mu S_\nu^*
&=&\delta_{d(\mu)-d(\nu),n}
p_vS_\mu S_\nu^*\nno
&=&p_v\Phi_nS_\mu S_\nu^*.\eean

As $\mu,\nu$ were arbitrary paths, this shows that $p_v\Phi_n$ is a
finite rank endomorphism. For arbitrary $a=\sum
c_jS_{\mu_j}S_{\nu_j}^*$, where the sum is finite, we may apply the same
reasoning to each $p_{s(\nu_j)}$ to see that $a\Phi_n$ is finite rank
for all $a\in A_c$.

To see that $a\Phi_k$ is compact for all $a\in A$, recall that
every $a\in A$ is a norm limit of a sequence $\{a_i\}_{i\geq
0}\subset A_c$. Thus for any $n\in\Z^k$
$a\Phi_n=\lim_{i\to\infty}a_i\Phi_n$ and so is compact.
\end{proof}

\begin{lemma}\label{compactendo} Assume that the $k$-graph $(\Lambda,d)$ is
locally finite and has no sinks. For all $a\in A$,
$a(1+\D^2)^{-1/2}$ is a compact
 endomorphism of the $F$-module $X$.
\end{lemma}

\begin{proof}
First let $a=p_v$ for $v\in \Lambda^0$. Then the sum
$$ R_{v,N}:=p_v\sum_{|n|=0}^N\Phi_n(1+n^2)^{-1/2}$$
is finite rank, by Lemma \ref{finrank}. We will show that the
sequence $\{R_{v,N}\}_{N\geq 0}$  is convergent with respect to
the operator norm $\Vert\cdot\Vert_{End}$ of endomorphisms of $X$.
Indeed, assuming that $M>N$, \begin{align} \Vert
R_{v,N}-R_{v,M}\Vert_{End}&=\Vert
p_v\sum_{|n|=N+1}^{M}\Phi_n(1+n^2)^{-1/2}\Vert_{End}\nno
&\leq(1+(N+1)^2)^{-1/2}\to 0,\end{align} since the ranges of the
$p_v\Phi_n$ are orthogonal for different $n$. Thus, using the
argument from Lemma \ref{finrank}, $a(1+\D^2)^{-1/2}\in
End^0_F(X)$ for all $a\in A_c$. Letting $\{a_i\}$ be a Cauchy
sequence from $A_c$, we have
$$\Vert a_i(1+\D^2)^{-1/2}-a_j(1+\D^2)^{-1/2}\Vert_{End}\leq\Vert
a_i-a_j\Vert_{End}=\Vert a_i-a_j\Vert_A\to 0,$$ since
$\Vert(1+\D^2)^{-1/2}\Vert\leq 1$. Thus the sequence
$a_i(1+\D^2)^{-1/2}$ is Cauchy in norm and
 we
see that $a(1+\D^2)^{-1/2}$ is compact for all $a\in A$.
\end{proof}

\begin{prop}\label{Kasmodule} Assume that the $k$-graph $(\Lambda,d)$ is
locally finite and has no sinks. Let $V=\D(1+\D^2)^{-1/2}$. Then
$(X,V)$ defines a class in $KK^{k\bmod2}(A,F)$.
\end{prop}

\begin{proof} We refer to \cite{K} for more information.
We need to show that various operators belong to $End^0_F(X)$.
First, $V-V^*=0$, so $a(V-V^*)$ is compact for all $a\in A$. Also
$a(1-V^2)=a(1+\D^2)^{-1}$ which is compact from Lemma
\ref{compactendo} and the boundedness of $(1+\D^2)^{-1/2}$.
Finally, we need to show that $[V,a]$ is compact for all $a\in A$.
First we suppose that $a\in A_c$. Then we have
 \bean
[V,a]&=&
[\D,a](1+\D^2)^{-1/2}-\D(1+\D^2)^{-1/2}[(1+\D^2)^{1/2},a](1+\D^2)^{-1/2}\nno
&=&b_1(1+\D^2)^{-1/2}+Vb_2(1+\D^2)^{-1/2},\eean where
$b_1=[\D,a]\in A_c$ and $b_2=[(1+\D^2)^{1/2},a]$. Provided that
$b_2(1+\D^2)^{-1/2}$ is a compact endomorphism, Lemma
\ref{compactendo} will show
that $[V,a]$ is compact for all $a\in A_c$. So consider the action
of $[(1+\D^2)^{1/2},S_\mu S_\nu^*](1+\D^2)^{-1/2}$ on
$x=\sum_{n\in\Z^k}x_n$. We
find \bea &&
\sum_{n\in\Z^k}[(1+\D^2)^{1/2},S_\mu S_\nu^*](1+\D^2)^{-1/2}x_n\nno
&=&\sum_{n\in\Z^k}
\left((1+(d(\mu)-d(\nu)+n)^2)^{1/2}-(1+n^2)^{1/2}\right)(1+n^2)^{-1/2}S_\mu
S_\nu^*x_n\nno &=&\sum_{n\in\Z^k}f_{\mu,\nu}(n)S_\mu S_\nu^*
\Phi_nx.\label{limit}\eea The function \ben
f_{\mu,\nu}(n)=
\left((1+(d(\mu)-d(\nu)+n)^2)^{1/2}-(1+n^2)^{1/2}\right)(1+n^2)^{-1/2}\een
goes to zero as $n^2\to\infty$. As the $S_\mu S_\nu^*\Phi_n$ are
finite rank with orthogonal ranges (for different $n$),
the sum in (\ref{limit})
converges in the endomorphism norm, and so converges to a compact
endomorphism. For general $a\in A_c$ we write $a$ as a finite
linear combination of generators $S_\mu S_\nu^*$, and apply the
above reasoning to each term in the sum to find that
$[(1+\D^2)^{1/2},a]$ is a compact endomorphism for all $a\in A_c$.

Now let $a\in A$ be the norm limit of a Cauchy sequence
$\{a_i\}_{i\geq 0}\subset A_c$. Then
$$\Vert[V,a_i-a_j]\Vert_{End}\leq 2\Vert a_i-a_j\Vert_{End}\to 0,$$
so the sequence $[V,a_i]$ is also Cauchy in norm, and so the limit
is compact.

It is also clear from the construction that if $k$ is even, the
Kasparov module is even (with grading given by $\omega_\C$)
and so belongs to $KK^0(A,F)$, while when $k$ is odd, the Kasparov
module belongs to $KK^1(A,F)$.
\end{proof}

\section{The Gauge Spectral Triple of a $k$-Graph
Algebra}\label{secondconstruction}

In this section we will construct a semifinite spectral triple for
those locally convex $k$-graph $C^*$-algebras which possess a
faithful, semifinite, lower-semicontinuous,
gauge invariant trace, $\tau$. Recall from
Proposition~\ref{trace=graphtrace} that such traces arise from
 faithful $k$-graph traces.

We will begin with the right $F_c$ module $X_c$. In order to deal
with the spectral projections of $\D$ we will also assume
throughout this section that $(\Lambda,d)$ is  locally finite and
has no sinks. This ensures, by Lemma \ref{finrank} that for all
$a\in A$ and $n\in\Z^k$ the endomorphisms $a\Phi_n$ of $X$ are
compact endomorphisms.

We define a ${\C}$-valued inner product on $X_c$ by
$$ \la x,y\ra:=\tau((x|y)_R)=\sum_{j=1}^{2^{[k/2]}}\tau(\Phi(x_j^*y_j))
=\sum_{j=1}^{2^{[k/2]}}\tau(x_j^*y_j).$$ Observe that this inner
product is linear in the second variable. We define the Hilbert
space $\HH=L^2(X,\tau)$ to be the completion of $X_c$ in the norm
coming from the inner product.

\begin{lemma}\label{endoproof} The $C^*$-algebra $A=C^*(\Lambda)$ acts on $\HH$
by an extension of left multiplication. This defines a faithful nondegenerate
$*$-representation of $A$. Moreover, any endomorphism of $X$ leaving $X_c$
invariant extends uniquely to a bounded linear operator on $\HH$.
\end{lemma}

\begin{proof} The first statement follows from the proof of
Proposition \ref{trace=graphtrace}. Now let $T$
be an endomorphism of  $X$ leaving $X_c$ invariant. Then
\cite[Cor 2.22]{RW},
$$(Tx|Ty)_R\leq \| T\|_{End}^2(x|y)_R$$
in the algebra $F$. Now the  norm of $T$ as an operator on $\HH$,
denoted $\Vert T\Vert_\infty$, can be computed in terms of the
endomorphism norm of $T$ by \begin{align}
\|T\|_\infty^2&:=\sup_{\|x\|_\HH\leq 1}\la
Tx,Tx\ra=\sup_{\|x\|_\HH\leq 1}\tau((Tx|Tx)_R)\nno &\leq
\sup_{\|x\|_\HH\leq 1}\| T\|_{End}^2\tau((x|x)_R)=\|
T\|_{End}^2.\end{align}
\end{proof}

\begin{cor} The endomorphisms $\{\Phi_n\}_{n\in\Z^k}$ define
mutually orthogonal projections on $\HH$. For any $K\subset \Z^k$
the sum $\sum_{n\in K}\Phi_n$ converges strongly to a projection
$\Phi_K$ in $\B(\HH)$. The projection $\Phi_{\ZZ^k}$ corresponding to
$K = \ZZ^k$ is equal to $\id_\HH$, so that for all $x\in \HH$ the
sum $\sum_n\Phi_nx$ converges in norm to $x$.
\end{cor}
\begin{proof} As in Lemma \ref{phiendo}, we can use the
continuity of the $\Phi_n$ on $\HH$, which follows from Lemma
\ref{endoproof}, to see that the relation
$\Phi_n\Phi_m=\delta_{n,m}\Phi_n$ extends from $X_c\subset\HH$ to
$\HH$. The proof of the strong convergence of sums of $\Phi_n$'s
is just as in Lemma \ref{phiendo} after  replacing the
$C^*$-module norm with the Hilbert space norm.
\end{proof}
\begin{lemma} The operator $\D$ extends to a closed unbounded self-adjoint
operator on $\HH$. The closure of the operator $\D|_{X_c}$ is
$\D$.
\end{lemma}
\begin{proof} The proof is essentially the same as the
$C^*$-module version, Lemma \ref{CstarDee}. By replacing the
$C^*$-module norm and the $C^*$-Cauchy-Schwartz inequality with
the Hilbert space analogues, the proof that $\D$ is closed goes
through as before. We then define $\mbox{dom}\ \D$ to be the
completion of $X_c$ in the norm
$$x\to\Vert x\Vert_{\HH,\D}:=\Vert
x\Vert_{\HH}+\Vert\D x\Vert_{\HH}.$$ The proofs of symmetry and
self-adjointness now follow just as in the $C^*$-module case. The
last statement follows from the definition of $\mbox{dom}\ \D$.
\end{proof}

The Hilbert space $\HH$ and operator $\D$ are two of the
ingredients of our spectral triple. We also need a $*$-algebra. In
fact $A_c$ will do the job, but it also has a natural completion
$\A$  which is useful too. To prove both these assertions we need
the following lemma. The proof is the same as \cite[Lemma 5.4]{PRen}.

\begin{lemma}\label{deltacomms} Let $\HH,\D$ be as above and let
$\dd=\sqrt{\D^*\D}=\sqrt{\D^2}$ be the absolute value of $\D$.
Then for $S_\al S_\beta^*\in A_c$, the operator $[\dd,S_\al
S_\beta^*]$ is well-defined on $X_c$, and extends to a bounded
operator on $\HH$ with
$$\Vert[\dd,S_\al S_\beta^*]\Vert_{\HH}\leq
\Bigl|d(\al)-d(\beta)\Bigr|.$$ Similarly, $\Vert[\D,S_\al
S_\beta^*]\Vert_\HH= \Bigl|d(\al)-(\beta)\Bigr|$.
\end{lemma}

\begin{cor}\label{smodense} The algebra $A_c$ is contained in the smooth domain
of the derivation $\delta$ where for $T\in\B(\HH)$,
$\delta(T)=[\dd,T]$. That is
$$ A_c\subseteq\bigcap_{n\geq 0}{\rm dom}\ \delta^n.$$
\end{cor}
\begin{defn} Define the $*$-algebra $\A\subset A$ to be the
completion of $A_c$ in the $\delta$-topology. By Lemma \ref{smo},
$\A$ is Fr\'{e}chet and stable under the holomorphic functional
calculus.
\end{defn}
\begin{lemma}\label{smoalg} If $a\in\A$ then $[\D,a]\in\A$ and the operators $\delta^k(a)$,
$\delta^k([\D,a])$ are bounded for all $k\geq 0$. If $\phi\in\A$
satisfies $\phi a=a=a\phi$, then $\phi[\D,a]=[\D,a]=[\D,a]\phi$.
The norm closed algebra generated by $\A$ and $[\D,\A]$ is $A\otimes
M_{2^{[k/2]}}(\C)$. In
particular, $\A$ is quasi-local.
\end{lemma}
We leave the straightforward proofs of these statements to the
reader.

At this point we have most of the structure required to define a
semifinite local spectral triple. The one remaining piece of
information we require is the compactness of $a(\lambda-\D)^{-1}$,
$\lambda\in\C\setminus\R$, $a\in\A$, relative to some trace on
some von Neumann algebra to which $\D$ is affiliated. There is a
canonical choice of von Neumann algebra and trace, and for this
choice $a(1+\D^2)^{-k/2}$ is in the domain of the Dixmier trace
for all $a\in\A$.

\subsection{Traces and Compactness Criteria}
We continue to assume that $(\Lambda,d)$ is a locally convex
locally finite $k$-graph  with no sinks and that $\tau$ is a
faithful, semifinite, lower-semicontinuous,
gauge invariant trace on $C^*(\Lambda)$.
 We will define a von Neumann algebra $\cn$
with a faithful semifinite normal trace $\tilde\tau$ so that
$\A\subset\cn\subset\B(\HH)$, where $\A$ and $\HH$ are as defined in the last
subsection. Moreover the operator $\D$ will be affiliated to $\cn$.
To state the theorem, we need some preliminary definitions and results.
\begin{defn}\label{v-N alg} Let $End^{00}_F(X_c)$ denote the algebra of finite rank
operators on $X_c$ acting on $\HH$. Define
$\cn=(End^{00}_F(X_c))''$, and let $\cn_+$ denote the positive
cone in $\cn$.
\end{defn}

\begin{defn}\label{t-def} Let  $T\in\cn$. For $n \in \N^k$ and $v \in \Lambda^0$,
let $|v|_n$ denote the number of paths of degree $n$ with source $v$.
Let $\Lambda \times_s^{\rm min} \Lambda$ denote the set of pairs
\[
\{(\alpha,\beta) \in \Lambda : s(\alpha) = s(\beta),
d(\alpha) \wedge d(\beta) = 0\}.
\]
For $(\alpha,\beta) \in \Lambda \times_s^{\rm min} \Lambda$, define
\[
\omega_{\alpha,\beta} (T) =
\frac{1}{|s(\alpha)|_{d(\beta)}}
\la s_\alpha s^*_\beta, T s_\alpha s^*_\beta\ra.
\]
Note that if $d(\alpha) = d(\beta) = 0$, then $\alpha = \beta = v$
for some $v \in \Lambda^0$, and since $s_v = p_v$ by convention,
we have $\omega_{v,v}(T)=\la p_v,Tp_v\ra.$ Define
\begin{equation}\label{eq:tilde tau def}
\tilde\tau:\cn_+\to[0,\infty],\quad\text{by}\quad
\tilde\tau(T)=
\lim_{L\nearrow \Lambda \times_s^{\rm min} \Lambda}
\sum_{(\alpha,\beta) \in L}\omega_{\alpha,\beta}(T)
\end{equation}
where $L$ increases over the net of finite subsets of
$\Lambda \times_s^{\rm min} \Lambda$.
\end{defn}

{\bf Remarks}\begin{itemize}
\item[(1)] For $T,S\in\cn_+$ and $\lambda\geq 0$ we have
$$\tilde\tau(T+S)=\tilde\tau(T)+\tilde\tau(S)\ \ \
\mbox{and}\ \ \ \tilde\tau(\lambda T)=\lambda\tilde\tau(T)\ \ \mbox{where}\ \
0\times\infty=0.$$
\item[(2)] Note that for $\mu \in \Lambda$, we have
$\omega_{\mu,s(\mu)}(T) = \langle s_\mu, Ts_\mu\rangle$ and
$\omega_{s(\mu),\mu}(T) = \frac{1}{|s(\mu)|_{d(\mu)}} \langle s^*_\mu, T s_\mu^*\rangle$.
Consequently, if $\Lambda$ is a $1$-graph then for $\mu \in \Lambda \setminus \Lambda^0$,
the map $\omega_\mu$ of \cite[Definition~5.10]{PRen} is precisely $\omega_{\mu,s(\mu)} +
\omega_{s(\mu),\mu}$, while for $v \in \Lambda^0$, $\omega_v = \omega_{v,v}$.
In particular, for a $1$-graph, \eqref{eq:tilde tau def} is just a slightly more efficient
expression for the definition of $\tilde\tau$ of \cite[Definition~5.10]{PRen}.
\end{itemize}

\begin{thm}\label{mainthm} Let $(\Lambda,d)$ be a locally convex locally finite $k$-graph with no sinks, and let $\tau$ be a faithful semifinite trace on
$C^*(\Lambda)$. Let $\cn$ be as in Definition~\ref{v-N alg} and let $\tilde\tau : \cn^+ \to [0,\infty]$ be as in Definition~\ref{t-def}. Then
\begin{itemize}
\item[(1)] $\tilde\tau$ defines a faithful normal semifinite trace on $\cn$.
\item[(2)] $(\A,\HH,\D)$ is a $QC^\infty$ $(k,\infty)$-summable odd local semifinite spectral triple relative to $(\cn,\tilde\tau)$.
\item[(3)] For all $a\in \A$, the operator $a(1+\D^2)^{-1/2}$ is not trace class.
\end{itemize}
Suppose that $v \in \Lambda^0$ satisfies $v\Lambda^{\le n}
= v\Lambda^n$ for all $n \in \N^k$. Then
$$\tilde\tau_\omega(p_v(1+\D^2)^{-k/2})=C_k\tau(p_v),$$
where $\tilde\tau_\omega$ is any Dixmier trace associated to
$\tilde\tau$, and $$C_k=\frac{2^{[k/2]}vol(S^{k-1})}{k}$$
\end{thm}

{\bf Remark}
The hypothesis that $v\Lambda^{\le n} = v\Lambda^n$ for all $n \in \N^k$
is perhaps somewhat opaque. This theorem generalises \cite[Theorem~5.8]{PRen},
which requires that the vertex $v$ has ``no sinks downstream'' to ensure
that $p_v = \sum_{s(\alpha) = v, |\alpha| = n} s_\alpha s^*_\alpha$ for all
$n \in \N$. The hypothesis that $v\Lambda^{\le n} = v\Lambda^n$ for all
$n \in \NN^k$ has precisely the same effect (consider relation~(CK4)). Indeed
this is precisely the notion that the $\Lambda^{\le n}$ notation was developed
to capture: $\Lambda^{\le n}$ is supposed to consist of all paths of degree $n$
together with all paths whose degree is less than $n$ because they originate
at a source in direction $n$ \cite{RSY1}.

\begin{prop}\label{tildetau}
 The function $\tilde\tau:\cn_+\to[0,\infty]$ defines a faithful normal
semifinite trace on $\cn$. Moreover,
$$End_F^{00}(X_c)\subset\cn_{\tilde\tau}:=
{\rm span}\{T\in\cn_+:\tilde\tau(T)<\infty\},$$
the domain of definition of $\tilde\tau$, and
$$\tilde\tau(\Theta^R_{x,y})=\la y,x\ra=\tau(y^*x),\ \ \ x,y\in X_c.$$
\end{prop}

The proof of this important, but technical, result is extremely similar
to that of \cite[Proposition 5.11]{PRen}, differing only in the details of
the calculations establishing the analogue of \cite[Equation~(18)]{PRen}
and showing that $\tilde\tau(\Theta^R_{x,y}) = \tau(y^*x)$ for all $x,y$.

\begin{lemma}\label{tracepvphi} Let $(\Lambda,d)$ be a
locally convex locally finite $k$-graph with no sinks and a
faithful gauge invariant trace $\tau$ on $C^*(\Lambda)$. Let $v\in
\Lambda^0$ and $n\in\Z^k$. Then
$$ \tilde\tau(p_v\Phi_n)\leq \tau(p_v)$$
with equality when $v \Lambda^{\le p} = v\Lambda^p$.
\end{lemma}

\begin{proof}  Let $n = n_+ + n_-$ where $n_+ \geq 0$, $n_- \le 0$,
and $n_+ \vee -n_- = n_+-n_-$. By Lemma \ref{finrank} and
Proposition \ref{tildetau} we have
\bean \tilde\tau\left(p_v\Phi_n\right)&=&
\tilde\tau\left(p_v\sum_{d(\al)=n_+,d(\beta)=-n_-}\frac{1}{|s(\beta)|_{-n_-}}
\Theta_{S_\al S_\beta^*,S_\al S_\beta^*}\right)\nno
&=&\tau\left(\sum_{d(\al)=n_+,d(\beta)=-n_-}
\frac{1}{|s(\beta)|_{-n_-}}(S_\al S_\beta^*|p_vS_\al
S_\beta^*)_R\right)\nno
&=&
\tau\left(\sum_{d(\al)=n_+,d(\beta)=-n_-}\frac{1}{|s(\beta)|_{-n_-}}
\Phi(S_\beta S_\al^*p_vS_\al S_\beta^*)\right)\nno&=&
\tau\left(\sum_{d(\al)=n_+,d(\beta)=-n_-,r(\al)=v}
\frac{1}{|s(\beta)|_{-n_-}}S_\al
S_\beta^*S_\beta S_\al^*p_v\right)\nno&=&
\tau\left(\sum_{d(\al)=n_+,r(\al)=v}S_\al S_\al^*p_v\right).\eean
If there are no sources within $|n_+|$ of $v$, then
$\sum_{d(\al)=n_+,r(\al)=v}S_\al S_\al^*=p_{r(\al)}=p_v$. Otherwise
the sum on the right is strictly less than $p_v$. So
$$\tilde\tau(p_v\Phi_n)\leq\tau(p_v)$$
with equality when there are no sources within $|n_+|$ of $v$.
\end{proof}

\begin{prop}\label{Dixytilde=tau} Assume that the locally convex
$k$-graph $(\Lambda,d)$  is locally finite, has no sinks and has a
faithful gauge invariant trace on $C^*(\Lambda)$. For all $a\in
\A_c$ the operator $a(1+\D^2)^{-k/2}$ is in the ideal
$\LL^{(1,\infty)}(\cn,\tilde\tau)$. When $v \in \Lambda^0$ satisfies
$v\Lambda^{\le n} = v\Lambda^n$ for all $n \in \N^k$, we have
$$\tilde\tau_\omega(p_v(1+\D^2)^{-k/2})=\frac{2^{[k/2]}vol(S^{k-1})}{k}
\tau(p_v).$$
\end{prop}

\begin{proof} It suffices to show this for a projection $a=p_v$ for
$v\in \Lambda^0$,
and extending to more general $a\in A_c$ using the arguments of
Lemma \ref{finrank}. We compute the partial sums defining the
trace of $p_v(1+\D^2)^{-k/2}$. Lemma \ref{tracepvphi} gives us
\begin{equation}
\tilde\tau\left(p_v\sum_{|n|\leq N}(1+n^2)^{-k/2}\Phi_n\right)
\leq\sum_{|n|\leq N}(1+n^2)^{-k/2}\tau(p_v).\label{dblestarry}\end{equation}
We have equality when $v\Lambda^{\le n} = v\Lambda^n$ whenever $|n| \le N$.
Since $\Lambda$ has no sinks, the sequence
$$\frac{1}{\log|\{n:|n|\leq N\}|}\sum_{|n|\leq N}(1+n^2)^{-k/2}
\tilde\tau(p_v\Phi_k)$$
is bounded (there is at least one `direction' in which $n$ can increase
indefinitely, so the sequence does not go to zero).  Hence $p_v(1+\D^2)^{-k/2}\in\LL^{(1,\infty)}$ and for
any $\omega$-limit we have
$$\tilde\tau_\omega(p_v(1+\D^2)^{-k/2})\leq
\omega\mbox{-}\!\lim\frac{2^{[k/2]}vol(S^{k-1})}{k\log
m}\sum_{m=0}^N(1+m^2)^{-1/2}\tilde\tau(p_v\Phi_k).$$ When
there are no sources in $\Lambda$, we have equality in Equation
(\ref{dblestarry}) for any $v\in \Lambda^0$ and so
$$\tilde\tau_\omega(p_v(1+\D^2)^{-k/2})=\frac{2^{[k/2]}vol(S^{k-1})}{k}
\tau(p_v).$$
\end{proof}

Computing the Dixmier trace when $v\Lambda^{\le n}$ may be strictly larger
than $v\Lambda^n$ for some $n$ is harder.

 {\bf Remark}  Using Proposition \ref{wellbehaved}, one can
check  that
 \be res_{s=0}\tilde\tau(p_v(1+\D^2)^{-k/2-s})=
\frac{k}{2}\tilde\tau_\omega(p_v(1+\D^2)^{-k/2}).\label{res}\ee We
will require this formula when we apply the local index theorem.

\begin{cor}\label{compactresolvent} Assume $\Lambda$ is locally finite,
has no sources and has a faithful $k$-graph trace. Then for all $a\in
A$, $a(1+\D^2)^{-1/2}\in\K_\cn$.
\end{cor}

\begin{proof} (of Theorem \ref{mainthm}.) That we
have a $QC^\infty$ spectral triple follows from Corollary
\ref{smodense}, Lemma \ref{smoalg} and Corollary
\ref{compactresolvent}. The properties of the von Neumann algebra
$\cn$ and the trace $\tilde\tau$ follow from Proposition \ref{tildetau}. The
$(k,\infty)$-summability and the value of the Dixmier trace comes
from Proposition \ref{Dixytilde=tau}. The locality of the spectral
triple follows from Lemma \ref{smoalg}.
\end{proof}

\section{The Local Index Theorem for the Gauge Spectral Triple}\label{index}

The local index theorem for semifinite spectral triples described
in \cite{CPRS2,CPRS3} is relatively simple for the spectral
triples constructed here. This is because of the simple way in
which the triples are built using the Clifford algebra.

In the following discussion we assume we have a fixed locally finite
locally convex $k$-graph $(\Lambda,d)$ without sinks, and possessing a
faithful $k$-graph trace. We let
$(\A,\HH,\D)$ be the associated gauge spectral triple relative to
$(\cn,\tilde\tau)$  constructed in
the previous section.

Elementary manipulations with the Clifford variables, like those
in \cite[Section 11.1]{BCPRSW}, along with the Dixmier trace results,
show that when $k$ is odd we are left with only one term in the
local index theorem
$$ \phi_k(a_0,a_1,...,a_k)=-\sqrt{2\pi i}\frac{1}{k!}\frac{1}{\sqrt{\pi}}
\Gamma(k/2+1)res_{r=(1-k)/2}
\tilde\tau(a_0[\D,a_1]\cdots[\D,a_k](1+\D^2)^{-(k-1)/2-r}).$$

When $k$ is even we are left with only two terms:
$$ \phi_k(a_0,a_1,...,a_k)=\frac{1}{k!}\Gamma(k/2+1)res_{r=(1-k)/2}
\tilde\tau(\gamma a_0[\D,a_1]\cdots[\D,a_k](1+\D^2)^{-(k-1)/2-r}),$$
$$\phi_0(a_0)=res_{r=(1-k)/2}\frac{1}{(r-(1-k)/2)}\tilde\tau(\gamma
a_0(1+\D^2)^{-(k-1)/2-r}).$$

The zero component in the even case likewise vanishes for our
examples. The reason for this is simply that we have complete symmetry
between the $\pm1$ eigenspaces of $\gamma$, and so for $Re(r)$ large
$$\frac{1}{(r-(1-k)/2)}\tilde\tau(\gamma a_0(1+\D^2)^{-(k-1)/2-r})=0.$$
Hence in this particular case, the local index theorem is in fact
computed using the Hochschild class (top component) of the Chern character,
\cite{CPRS1}.

In Proposition~\ref{prp:sufficient}, we will describe a class of
$k$-graphs which admit faithful graph traces. For full details, see
Appendix~\ref{appendix}; for the time being we need only two facts
established there: (1) that for such $k$-graph, the $K$-theory of
$C^*(\Lambda)$ resides entirely on the set of ideals of
$C^*(\Lambda)$ corresponding to ends (see Definition~\ref{dfn:ends})
of $\Lambda$; and (2) that for each end of $\Lambda$, the associated
ideal is of the form $\K \otimes C(\T^l)$ for some $0 \le l \le k$.
In particular~(1) implies that it is only necessary to produce
generators of $K$-theory corresponding to these ends.

The form of the Chern character given above shows that in odd
dimensions we can detect only ends for which the number $l$ in~(2)
above is odd, whilst in even dimensions we can only detect ends for
which $l$ is even. A simple analysis based on the Clifford algebra
then shows that in fact we can only pair with ends where $l = k$;
that is, ends which are $k$-tori, $k$-planes, or more generally
$k$-cylinders.

Before producing an example of what this kind of index pairing can
tell us, we discuss the relationship between the $KK$-index
pairing with values in $K_0(F)$ and the semifinite index theorem.

\begin{thm}\label{thm:kk=res} Let $\Lambda$ be a locally convex, locally finite
$k$-graph without sinks which admits a faithful graph trace, let
$\tau$ be the corresponding semifinite trace on $A = C^*(\Lambda)$,
and let $(\A,\HH,\D)$ be the gauge spectral triple $($relative to
$(\cn, \tilde\tau))$ obtained from Theorem~\ref{mainthm}. Let
$(X,\D)$ be the corresponding Kasparov module with class in
$KK^k(A,F)$. Let $x\in K_k(A)$ be a $K$-theory class. Then
$$\tilde\tau_*([x\times(X,\D)])=Ch_{(\A,\HH,\D)}(Ch(x)).$$
\end{thm}

\begin{proof}
Let us first consider the even case, where we have the $K_0(F)$-valued
index of $p\D_+p$ on $X$, where $p$ is a projection in $A$.
The projections defining $\ker(p\D_+p)$ and $\mbox{coker}(p\D_+p)$
are compact endomorphisms of the module $X$, and moreover map $X_c$ to
itself. This last assertion follows because $\D$ maps $X_c$ to itself,
and $p$ may be chosen to lie in $A_c$ which preserves $X_c$. The
reason we can do this is that $K_0(A)=\lim K_0(\phi_nA\phi_n)$ where
$\phi_n$ is any local approximate unit for $A_c$, \cite{R1}. Hence the
kernel and cokernel projections are actually endomorphisms preserving
$X_c$.

Now such endomorphisms extend to act on the Hilbert space in a unique
way. Since $\HH=\overline{X_c}$ with respect to the norm coming from
the inner product, we see that the Hilbert space kernel and cokernel
projections are given by the extension of the $C^*$-module
projections. So we have, by Lemma \ref{tracepvphi},
$$\tilde\tau-\mbox{Index}(p\D_+p)=
\tilde\tau(Q_{\ker(p\D_+p)}-Q_{{\rm coker}(p\D_+p)})=
\tilde\tau_*([\mbox{Index}(p\D_+p)])=\tilde\tau_*([p]\times[(X,\D)]).$$

The argument for the odd pairing is now exactly the same, except that
we consider the kernel and cokernel projections of $PuP$ where $P$ is
the non-negative spectral projection of $\D$ and $u$ is  unitary.
The upshot is that
$$
\tilde\tau-\mbox{Index}(PuP)=\tilde\tau_*([\mbox{Index}(PuP)])=\tilde\tau_*([u]\times[(X,\D)]).$$

Now we wish to relate the $\tilde\tau$ index to the pairing of Chern
characters. However, by \cite{CPRS4}, this is precisely the main theorems of
\cite{CPRS2,CPRS3} in the odd and even cases respectively, and so we
are done.
\end{proof}

We will conclude with an example which indicates the kinds of
information one might hope to obtain from the semifinite index
theorem. In order to present the example explicitly, we first
produce representatives for generators of $K$-theory
coming from ends of graphs satisfying Proposition
\ref{prp:sufficient}. For this we need generators of the
$K$-theory of ordinary tori (those for planes are of course well known).
In fact we really only need those
generators which pair with the Dirac class. These in turn can all be
obtained, using the universal coefficient theorem and the fact that
the $K$-theory of tori is free abelian, by using iterated products
with the circle.

We illustrate this with a specific example; the nontrivial generator
of $K^2(\T^2)=\Z^2$ (the other products are simpler). We wish to
compute the product of $[u]\in K_1(C(\T^1))$ with itself to obtain the
nontrivial element of $K_0(\T^2)$. Let $\theta,\phi\in[0,2\pi]$,
and set $z=e^{i\theta}$. Then $\theta\to z$ represents the generator
$[u]$.
Define
$$ K(\theta)=\bma 0 & z\\ \bar z & 0\ema,\ \ \ S=\bma 0 & 1\\ 1 & 0\ema,\ \
\ Y(\phi,\theta)=e^{i\phi K(\theta)/4}e^{iS/4}.$$

Then the product $[u]\times[u]$ is the class of the projection
$$P(\phi,\theta)=Y(\phi,\theta)^*\bma 1 & 0\\ 0 & 0\ema
Y(\phi,\theta).$$
A lengthy computation shows that
$$P(\phi,\theta)=\bma 1-\sin^2(\phi/2)\cos^2(\theta/2)&
\frac{i}{2}\sin(\phi)\cos^2(\theta/2)-\frac{1}{2}\sin(\phi/2)\sin(\theta)\\
-\frac{i}{2}\sin(\phi)\cos^2(\theta/2)-\frac{1}{2}\sin(\phi/2)\sin(\theta)
& \sin^2(\phi/2)\cos^2(\theta/2).\ema$$
An even lengthier calculation using the residue cocycle from the local
index theorem, \cite{CPRS2,CPRS3}, shows that
$P(\phi,\theta)$ has pairing with
$$ \mbox{Dirac}_{\T^2}=\left(C^\infty(\T^2), L^2(\T^2)\otimes\C^2,\bma 0 &
-\p_\phi+i\p_\theta\\ \p_\phi+i\p_\theta\ema\right)$$
equal to one. Hence $p$ is the desired generator.

\begin{figure}[ht]
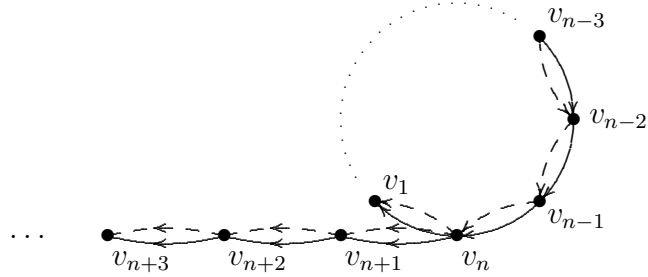

\[
\beginpicture
\setcoordinatesystem units <2.5em, 2.5em>
\put{$\bullet$} at 1.5 0
\put{$v_{n+3}$}[lt] at 1.55 -0.2
\put{$\bullet$} at 3 0
\put{$v_{n+2}$}[lt] at 3.05 -0.2
\put{$\bullet$} at 4.5 0
\put{$v_{n+1}$}[lt] at 4.55 -0.2
\put{$\bullet$} at 6 0
\put{$v_n$}[lt] at 6.05 -0.2
\put{$\bullet$} at 7.06 0.44
\put{$v_{n-1}$}[lt] at 7.16 0.34
\put{$\bullet$} at 7.5 1.5
\put{$v_{n-2}$}[l] at 7.7 1.5
\put{$\bullet$} at 7.06 2.56
\put{$v_{n-3}$}[lb] at 7.16 2.66
\put{$\bullet$} at 4.94 0.44
\put{$v_{1}$}[lb] at 5.04 0.54
\put{ } at 0 3
\put{$\dots$} at 0.5 0
\setdots
\circulararc 150 degrees from 6.75 2.8 center at 6 1.5
\setsolid
\circulararc 135 degrees from 6 0 center at 6 1.5
\arr(6.08, 0.002)(6.07,0.001)
\arr(6.08, 0.071)(6.07,0.061)
\arr(7.12, 0.51)(7.11, 0.50)
\arr(7.066, 0.54)(7.065, 0.53)
\arr(7.49, 1.59)(7.49, 1.58)
\arr(7.43, 1.58)(7.44, 1.57)
\circulararc -45 degrees from 6 0 center at 6 1.5
\arr(5.01, 0.38)(5, 0.39)
\arr(5.02, 0.44)(5.01, 0.441)
\setquadratic
\setdashes
\plot 6 0 6.48 0.32 7.06 0.44 /
\plot 7.06 0.44 7.18 1.02 7.5 1.5 /
\plot 7.06 2.56 7.18 1.98 7.5 1.5 /
\plot 6 0 5.52 0.32 4.94 0.44 /
%\setsolid
\plot 1.5 0 2.25 0.1 3.0 0 /
\arr(2.101, 0.1)(2.10, 0.1)
\arr(2.101, -0.1)(2.10, -0.1)
\plot 3.0 0 3.75 0.1 4.5 0 /
\arr(3.601, 0.1)(3.60, 0.1)
\arr(3.601, -0.1)(3.60, -0.1)
\plot 4.5 0 5.25 0.1 6.0 0 /
\arr(5.101, 0.1)(5.10, 0.1)
\arr(5.101, -0.1)(5.10, -0.1)
\setsolid
\plot 1.5 0 2.25 -0.1 3.0 0 /
\plot 3.0 0 3.75 -0.1 4.5 0 /
\plot 4.5 0 5.25 -0.1 6.0 0 /
\endpicture
\]
\caption{The $2$-graph $\Lambda_n$}
\label{fig:2-torus}
\end{figure}

\textbf{Example.}
Consider the $2$-graph $\Lambda_n$ whose skeleton is illustrated in
Figure~\ref{fig:2-torus}

We label the solid edge whose range is $v_i$ by $e_i$ and the dashed
edge with the same range is labelled $f_i$ for all $i < n$. Without
the infinite tail to the left, we think of this $2$-graph as the
'$n$-point $2$-torus'. To justify this, note that $\Lambda_n$ is the
pull-back of the graph $E_n$ from \cite{PRen} with respect to the
functor $\lambda \mapsto d(\lambda)_1 + d(\lambda)_2$, and so
\cite[Corollary~3.5(iii)]{KP} shows that
\[
C^*(\Lambda_n) \cong C^*(E_n) \otimes C(\T) \cong \K \otimes C(\T)
\otimes C(\T) \cong \K \otimes C(\T^2);
\]
in particular $p_{v_n}$ is a full projection in $C^*(\Lambda_n)$ and
$p_{v_n} C^*(\Lambda_n) p_{v_n}$ is isomorphic to $C(\T^2)$.

We wish to describe the isomorphism explicitly. To do this, first notice that
$p_{v_n} C^*(\Lambda) p_{v_n}$ is generated by the unitaries
\[
u_1 := s_{e_1} s_{e_2} \dots s_{e_n}
\quad\text{and}\quad
u_2 = s_{f_1} s_{e_2} \dots s_{e_n}.
\footnote{Note that while $w = s_{f_1} s_{f_2} \cdots s_{f_n}$ would appear to be
a more natural candidate for the second generator than $u_2$, it is easy to see
that $u_2$ does not belong to $C^*(\{u_1,w\})$ and so $u_1$ and $w$ do not generate
$p_{v_n} C^*(\Lambda) p_{v_n}$. From a $K$-theoretic point of view, however, the
distinction is not important: one can check that $u_2$ and $w$ have the same class
in $K_1(p_{v_n} C^*(\Lambda) p_{v_n})$.}
\]
For $1 \le i < j$ let $\theta_{i,j} := s_{e_j} \dots s_{e_i+1}$, for
$i > j$,  let $\theta_{i,j} := \theta_{j,i}^* = s_{e_{j+1}}^* \dots s_{e_i}^*$, and
for $i = j$, let $\theta_{i,i} := p_{v_i}$. For $i = 1,2$ define
$U_i \in \Mm(C^*(\Lambda_n))$ by $U_i := \sum_{j \in \N} \theta_{n,j} u_i \theta_{j,n}$.

\begin{lemma}\label{lem:Lambda_n}
The $C^*$-algebra $C^*(\Lambda_n)$ is generated by the elements
$\{\theta_{i,j} : i,j \in \N\}$ together with $\{U_1, U_2\}$. The
$\theta_{i,j}$ form a system of nonzero matrix units, the unitaries
$U_1, U_2$ commute and have full spectrum, and hence there is an
isomorphism of $C^*(\Lambda_n)$ onto $\K \otimes C(\T^2)$ which
takes $\theta_{i,j}U_1^mU_2^n$ to the function $(w,z) \mapsto
\Theta_{i,j} \otimes w^m z^n$. Moreover, the core $F_n =
C^*(\Lambda_n)^\gamma$ is isomorphic to $\oplus^{n}_{j=1} \K$.
\end{lemma}
\begin{proof}
The Cuntz-Kreiger relations show that the $\theta_{i,j}$ are matrix
units and that $U_1$ and $U_2$ are unitaries. Since the
$\theta_{n,j}$ have orthogonal range projections, $C^*(\{U_1,
U_2\})$ is canonically isomorphic to $C^*(\{u_1, u_2\})$, which in
turn is canonically isomorphic to $C(\T^2)$ (see
Proposition~\ref{prp:bdry path algebra}). It is easy to check that
$U_1$ and $U_2$ commute with the matrix units $\theta_{i,j}$ so
$C^*(\{U_1, U_2, \theta_{i,j} : i,j \in \N\}) \cong \K \otimes
C(\T^2)$ (it is worth noting that compression by $p_{v_n} =
\theta_{n,n}$ takes $U_i$ to $u_i$). It now remains to show that
this algebra is all of $C^*(\Lambda_n)$.

For $i \not= 1$, we have $s_{e_i} = \theta_{i,i+1}$, and we have
$s_{e_1} = u_1\theta_{n,1}$ and $s_{f_2} = u_2\theta_{n,1}$. The only possible
factorisation rule for $\Lambda_n$ satisfies $e_{i+1}f_i = f_{i+1}e_i$ for all
$i$, and it now follows that $s_{f_i} = \theta_{1,i}s_{f_1}\theta_{i-1,n}$ for
all $i > 1$. Hence all the generators of $C^*(\Lambda_n)$ belong to
$C^*(\{U_1, U_2, \theta_{i,j} : i,j \in \N\})$ and it follows that $C^*(\Lambda_n)$
is isomorphic to $\K \otimes C(\T^2)$ as required.

To see that $F_n$ is isomorphic to $\oplus^{n}_{i=1} \K$, first
observe that the subalgebra
$C^*(\{s_\alpha : \alpha \in \Lambda_n : d(\alpha)_2 = 0\})$
generated by paths consisting only of solid edges is canonically
isomorphic to
the graph algebra $C^*(E_n)$ described in \cite[Corollary~6.6]{PRen},
and that this
isomorphism intertwines the restriction of the
gauge action on $C^*(\Lambda_n)$ to
$(\T,1)$ and the gauge action on $C^*(E_n)$.

It is shown in
\cite{PRen} that the core of $C^*(E_n)$ is isomorphic to $\oplus^{n}_{i=1} \K$:
the minimal projections in the $j^{\rm th}$ copy of $\K$ are the
vertex
projections
$s_{v_i} : i \cong l \mod n$, and for $i \ge j$, the $(i,j)^{\rm th}$ matrix unit is
$\theta^l_{i,j} := s_\eta s_{L(v_j)^{i-j}}^* s^*_\zeta$ where $\eta$ is the shortest path
from $v_l$ to $v_{in + l}$, $\zeta$ is the shortest path from $v_l$ to $v_{jn + l}$, and
$L(v_l)$ is the loop of length $n$ based at $v_l$.

Hence it suffices to show here that
\begin{equation}\label{eq:one dimension}
F_n =
\clsp\{s_\mu s^*_\nu : d(\mu) = d(\nu), d(\mu)_2 = d(\nu)_2 = 0, s(\mu) = s(\nu)\}.
\end{equation}
Recall from \cite[Section~4.1]{RSY1} that
\[
F_n = \clsp\{s_\mu s^*_\nu : \mu,\nu \in \Lambda_n, d(\mu) = d(\nu), s(\mu) = s(\nu)\},
\]
so we just need to show that if $\mu,\nu \in \Lambda_n$ satisfy $d(\mu) = d(\nu)$ and
$s(\mu) = s(\nu)$, then there exist $\eta$ and $\zeta$ such that $d(\eta) = d(\zeta) = (c,0)$
for some $c \in \N$, $s(\eta) = s(\zeta)$,  and $s_\eta s^*_\zeta = s_\mu s^*_\nu$.
Fix $\mu, \nu \in \Lambda_n$ with $d(\mu) = d(\nu)$ and $s(\mu) = s(\nu)$, and write
$(a,b)$ for $d(\mu)$. Let $\beta$ be the unique path of degree $(b,0)$ whose range is equal
to the source of $\mu$. By the factorisation property we can express $\mu = \mu_1\mu_2$ and
$\nu = \nu_1\nu_2$ where $d(\mu_1) = d(\nu_1) = (a,0)$ and $d(\mu_2) = d(\nu_2) = (0,b)$.
Applying the factorisation property again, we obtain
\[
\mu\beta = \mu_1\mu_2\beta = \mu_1\beta'\mu_2'
\quad\text{and}\quad
\nu\beta = \nu_1\nu_2\beta = \nu_1\beta''\nu_2'
\]
where $d(\beta') = d(\beta'') = (b,0)$ and $d(\mu_2') = d(\nu_2') = (0,b)$.
Since $|\mu_1\beta'| = |\mu_1| + |\beta| = a + b = |\mu|$, we have $s(\mu_1\beta') = s(\mu)$,
and similarly $s(\nu_1\beta'') = s(\nu) = s(\mu)$. Hence $\mu_2'$ and $\nu_2'$ are two
paths with the same degree and same range. Since $v\Lambda_n^p$ is a singleton for each $v$
and $p$, it follows that $\mu'_2 = \nu'_2$, so $s_{\mu'_2}s^*_{\nu'_2} = p_{s(\mu)}$ by~(CK4).
But now
\[
s_\mu s^*_\nu
= s_\mu p_{s(\mu)} s^*_\nu
= s_\mu s_\beta s^*_\beta s^*_\nu
= s_{\mu_1\beta'\mu'_2} s^*_{\nu_1\beta''\nu'_2}
= s_{\mu_1\beta'} s_{\mu'_2} s^*_{\mu'_2} s_{\nu_1\beta''}
= s_{\mu_1\beta'} s_{\nu_1\beta''}.
\]
Since $d(\mu_1\beta') = (a+b,0) = d(\nu_1\beta'')$ and since $s(\beta') = s(\beta'') = s(\mu)$,
this establishes~\eqref{eq:one dimension}.
\end{proof}

So for all $n$ we have $K_1(C^*(\Lambda_n))\cong
K_0(C^*(\Lambda_n))\cong\Z^2$. Choose any unitaries $v_1,v_2\in
(C^*(\Lambda_n)_c)^+$, the one-point unitization of the span of the
generators, such that $v_1,v_2$ represent the classes of the standard
generators $z_1,z_2$ of $K_1(C(\T^2))$. Then we obtain, as above, a
projection $P(v_1,v_2)$ representing the class of the Bott generator
in $ K_0(C^*(\Lambda_n))$. Using this, we may compute the pairing of
the Kasparov module $(X_n,\D_n)$ constructed for $C^*(\Lambda_n)$ with
the Bott projector.

As in \cite{PRen} we will compute first with the `$n$-point 2-torus',
the analogous calculation for the 2-graph $\Lambda_n$ will then follow
from the isomorphism $K_0(\K^n)\cong K_0(\C^n)=\Z^n$.

Let $\phi:C(\T^2)\to M_n(C(\T^2))$ be given by
$$\phi(f(z_1,z_2))=\theta_{n,n}f(w_1,w_2)\theta_{n,n}+(1-\theta_{n,n})=p_{v_n}f(w_1,w_2)p_{v_n}+(1-p_{v_n}).$$

Here we have set $w_1=u_1$ and $w_2=w$, as in the proof of Lemma
\ref{lem:Lambda_n}, and denoted the generating unitaries of
$C(\T^2)$ by $z_1,z_2$. Also $\theta_{n,n}$ is the projection
$p_{v_n}$. Let $(X,\D)$ be the Kasparov module for the $n$-point
2-torus built from the gauge action of $\T^2$. Then
$\D=\sum^n_{j=1}p_{v_j}\D=\sum p_{v_j}\D p_{v_j}$, and the pull-back
of $(X,\D)$ by $\phi$ is
$$\phi^*(X,\D)=(p_{v_n}X,p_{v_n}\D)\oplus\mbox{ degenerate module}\in
KK^0(C(\T^2),F)$$
since $1-p_{v_n}$ commutes with $\D$. The isomorphism $\psi:F\to\C^n$
given by
$$\psi(\sum_{j=1}^nz_jp_{v_j})=(z_1,\dots,z_n)$$
gives us
$$\psi_*\phi^*(X,\D)=\oplus_{j=1}^n(p_{v_n}Xp_{v_j},p_{v_n}\D)\in
\oplus_{j=1}^n K^0(C(\T^2).$$ The class of
$(p_{v_n}Xp_{v_j},p_{v_n}\D)$ is easily seen to be the Dirac
operator on $\T^2$ for the usual flat metric. In the following we
will identify $F$ with $\C^n$ (suppressing $\psi$).

Now we can
compute the pairing of $\D$ with $P(w_1,w_2)=\phi(P(z_1,z_2))$ where
$P(z_1,z_2)$ is the Bott projector of $\T^2$ constructed earlier. We
have
\begin{align} \la
[P(w_1,w_2)],[(X_n,\D_n)]\ra&=\la\phi_*([P(z_1,z_2)],[(X_n,\D_n)]\ra\nno
&=\la[P(z_1,z_2)],\phi^*[(X_n,\D_n)]\ra\ \ \mbox{functoriality of Kasparov
product}\nno
&=\la [P(z_1,z_2)],[(p_{v_n}X_n,p_{v_n}\D_n)]\oplus[\mbox{degenerate
module}]\ra\nno
&=\la P(z_1,z_2),\oplus_{j=1}^n[(p_{v_n}X_np_{v_j},p_{v_n}\D_n)]\ra
\nno
&=\oplus_{j=1}^n\la P(z_1,z_2),[\mbox{Dirac}_{\T^2}]\ra\nno
&=\oplus_{j=1}^n\la
[z_1]\times[z_2],[\mbox{Dirac}_{\T^1}]\times[\mbox{Dirac}_{\T^1}]\ra\
\ \mbox{by \cite[Theorem 10.8.7]{HR}}
\nno
&=-\oplus_{j=1}^n\la
[z_1],[\mbox{Dirac}_{\T^1}]\ra\la[z_2],[\mbox{Dirac}_{\T^1}]\ra\ \
\mbox{\cite[Chapter 9]{HR}}\nno
&=-(1,1,\dots,1)\in \Z^n=K_0(\C^n).\end{align}

Using Theorem \ref{thm:kk=res}, we may compute the pairing of the
Bott class with the spectral triple $(\A,\HH_n,\D_n)$ (where the
$k$-graph trace is chosen to be equal to $1$ on each vertex) by
applying $\tilde\tau_*$ to this last computation. We obtain
$$  \la
[P(w_1,w_2)],[(\A,\HH_n,\D_n)]\ra=-n.$$

The number $n$ appears basically because the multiplicity provided by
the core has given us $n$ copies of the Dirac operator at each point.

Now one can add the handle to the $n$-point 2-torus to get the 2-graph
$\Lambda_n$. The core becomes $\K^n$ and an argument entirely
analogous to the above shows again that the number $n$ emerges from the
pairing of $\mathcal{D}$ with the class of the Bott projector.

Note that this example can be generalised to an $n$-point $k$-torus with a ``handle". A
similar argument to that of Lemma~\ref{lem:Lambda_n} shows that the resulting $k$-graph
$\Lambda^k_n$ satisfies $C^*(\Lambda^k_n) \cong \K \otimes C(\T^k)$ independent of $n$, but
that the core $F^k_n$ is always isomorphic to $\oplus^{n}_{i=1} \K$. We can therefore see
that $n$ appears in the index computation in each case.

The point of this example is as follows. Whilst $C^*(\Lambda_n)\cong
C^*(\Lambda_m)$ for all $n,m$, as $k$-graph algebras they are,
somewhat vaguely, `different' for $n\neq m$. This difference is
embodied by equivalently the different gauge actions, the
nonisomorphic cores, and the different presentations as universal
algebras. Because $(X_n,\D_n)$ is constructed from the gauge action,
one would expect that $[(X_n,\D_n)]\in
KK^0(C^*(\Lambda_n),C^*(\Lambda_n)^\gamma)$ could `see' these
differences.

However, given a semifinite spectral triple $(\A,\HH,\D)$ relative
to $(\cn,\tau)$, we have no knowledge of the possible range of the
index; any real number is possible \emph{a priori}. What theorem
\ref{thm:kk=res} says, at least in this case, is that the semifinite
index {\em is} `quantised' --- the resulting index for
$C^*(\Lambda_n)$ is always a multiple of $n$. In \cite{KNR},
inspired by the result for $k$-graphs, it is shown that a similar
result is true for any semifinite spectral triple.

We list two results for this very simple example which indicate the
kinds of information one may draw from semifinite index theory in
general:

1) No combination of operator homotopy and addition of degenerate
spectral triples, \cite{CPRS1}, can make $(A_c,\HH_n,\D_n)$ and
$(A_c,\HH_m,\D_m)$ unitarily equivalent.
%\begin{proof} If we could make them unitarily equivalent by these
%procedures, the resulting index pairings would be equal, which they
%are not.
%\end{proof}

2) The gauge actions on $C(\T^2)\otimes \K$ coming from the
presentations as $C^*(\Lambda_n)$ and $C^*(\Lambda_m)$, $n\neq m$,
are not homotopic in $\Aut(C(T^2)\otimes K)$.
%\begin{proof} If they were homotopic, we could likewise deform the
%generators into each other, which by $1)$ we can not.
%\end{proof}

From the point of view of $k$-graph algebas, what is interesting
here is not the differences between $\Lambda_m$ and $\Lambda_n$, or
between the cores of the corresponding $C^*$-algebras, but rather
that the semifinite index can detect these differences in the
algebras. For example, while it is obvious that for $n\neq m$ there
is no gauge equivariant isomorphism $\phi:C^*(\Lambda_n)\to
C^*(\Lambda_m)$ (such a map would give an isomorphism $\K^n \cong
\K^m$ of the fixed point algebras), and hence there is no functorial
isomorphism between the $2$-graphs $\Lambda_n$ and $\Lambda_m$ (such
an equivalence would give rise to a gauge equivariant isomorphism of
$C^*$-algebras), it is not obvious that the semifinite index (which
sees the $C^*$-algebra and not the graph) should detect such
information.

Thus the semifinite index reflects finer details of the
$C^*$-algebra than the ordinary Fredholm index possibly could.
% This
%example is of course transparent; we have not learnt anything
%genuinely new about $C^*(\Lambda_n)$. However the example is
%illustrative of some possible uses of semifinite index theory.

\appendix

\section{$k$-graphs which admit faithful graph traces}\label{appendix}

In this appendix we formulate two necessary conditions
(Lemma~\ref{infpaths} and Corollary~\ref{infpaths to end}), and one
sufficient condition (Proposition~\ref{prp:sufficient}), for a
$k$-graph $\Lambda$ to admit a faithful graph trace. Our sufficient
condition is certainly much stronger than need be; indeed, the
$C^*$-algebra of a $k$-graph satisfying our condition is Morita
equivalent to a commutative $C^*$-algebra whereas, for example,
\cite{PRRS} contains many examples of $2$-graphs which admit
faithful graph traces and whose $C^*$-algebras are simple A$\T$
algebras with real rank 0. However, our condition is a direct
generalisation of the corresponding result in \cite{PRen} which has
already attracted independent interest. Moreover, the results about
$k$-graphs satisfying this condition, including the $K$-theory
calculations, are new and should be of independent interest to the
$k$-graph community.

For the purposes of our first two results we say that paths $\mu$
and $\nu$ in a $k$-graph $\Lambda$ are \emph{orthogonal} if the
range projections $s_\mu s^*_\mu$ and $s_\nu s^*_\nu$ are orthogonal
in $C^*(\Lambda)$. By \cite[Proposition~3.5]{RSY1}, $\mu$ and $\nu$
are orthogonal if and only if they have no common extensions.

\begin{lemma}[{cf \cite[Lemma~3.5]{PRen}}]\label{infpaths}
Suppose that $(\Lambda,d)$ is a row-finite $k$-graph and
there are vertices $v, w\in \Lambda^0$ with an infinite number of
mutually orthogonal paths from $w$ to $v$. Then there is no faithful
$k$-graph trace on $\Lambda^0$.
\end{lemma}
\begin{proof}
Let $(\lambda_n)_{n \in \N}$ be the infinite set of orthogonal paths
from $w$ to $v$. Suppose that $\tau$ is a trace on $C^*(\Lambda)$.
For each $n$, $\tau(s_{\lambda_n} s^*_{\lambda_n}) =
\tau(s^*_{\lambda_n} s_{\lambda_n}) = \tau(p_w)$. It follows that for
any $N$, we have $\tau(p_v) \ge \sum_{n=1}^N \tau(s_{\lambda_n}
s^*_{\lambda_n}) = N\tau(p_w)$, and it follows that $\tau(p_w) = 0$.
Hence $g_\tau(w) = 0$, and it follows from
Proposition~\ref{trace=graphtrace} that no $k$-graph trace on
$\Lambda^0$ is faithful.
\end{proof}

\begin{cor}[{cf \cite[Corollary~3.7]{PRen}}]\label{infpaths to end}
Suppose that $(\Lambda,d)$ is a row-finite $k$-graph and
there exists a vertex $v\in \Lambda^0$ with an infinite number of
mutually orthogonal paths from an end to $v$. Then there is no faithful
$k$-graph trace on $\Lambda^0$.
\end{cor}
\begin{proof}
Since $k$-graph traces are constant on ends by Remarks~\ref{rmk:const on ends},
the proof is identical to that of Lemma~\ref{infpaths}.
\end{proof}

We now aim to provide a sufficient condition for a $k$-graph to
admit a faithful $k$-graph trace.

\begin{not*}
Let $\Lambda$ be a locally convex row-finite $k$-graph. For ends $x$
and $y$ of $\Lambda$, we write $x \sim y$ if and only if $x(n) =
y(m)$ for some $n \le d(x)$ and $m \le d(y)$. This defines an
equivalence relation on ends of $\Lambda$, and we write $[x]$ for
the equivalence class of an end $x$ under $\sim$.

If a vertex $v$ lies on an end of $\Lambda$, then $v\Lambda^{\le
\infty} = \{x_v\}$, where $x_v$ is itself an end of $\Lambda$.
\end{not*}

\begin{prop}[{cf \cite[Propositions 3.8~and~3.9]{PRen}}]\label{prp:sufficient}
Let $\Lambda$ be a locally convex row-finite $k$-graph. Suppose that
there is a function $v \mapsto n_v$ from $\Lambda^0$ to $\N^k$ such
that for each $v \in \Lambda^0$ and each $\lambda \in v\Lambda^{\le n_v}$,
$s(\lambda)$ lies on an end of $\Lambda$.
\begin{itemize}
\item[(a)] If $g : \Lambda^0 \to \R^+$ is a $k$-graph trace,
then there is a well-defined function from
$\Ends(\Lambda)/\!\sim$ to $\R^+$ satisfying $g([x]) := g(x(0))$,
and
\begin{equation}\label{eq:graph trace formula}
g(v) = \sum_{\lambda \in v\Lambda^{\le n_v}} g([x_{s(\lambda)}])
\quad\text{ for every $v \in \Lambda^0$}.
\end{equation}
\item[(b)] Conversely, given any function
$g$ from $\Ends(\Lambda)/\!\sim$ to $\R^+$, there is a unique graph-trace
$\overline{g}$ on $\Lambda$ satisfying $\overline{g}(x(0)) = g([x])$
for all $x \in \Ends(\Lambda)$.
\end{itemize}
\end{prop}

Before proving the Proposition we need to know that for a fixed
function $g$ from $\Ends(\Lambda)/\!\sim$ to $\R^+$, the
formula~\eqref{eq:graph trace formula} is independent of the choice
of function $v \mapsto n_v$.

\begin{lemma} \label{lem:independent of n_v}
Suppose that $\Lambda$ satisfies the hypotheses of
Propopoition~\ref{prp:sufficient}, and let $g$ be a function from
$\Ends(\Lambda)/\!\sim$ to $\R^+$. Define $g(v) := g([x_v])$
for each vertex $v$ that lies on an end of $\Lambda$. Fix
$v \in \Lambda^0$ and suppose $n_1, n_2 \in \N^k$ each have the property
that $s(\lambda)$ lies on an end of $\Lambda$ for each
$\lambda \in v\Lambda^{\le n_i}$. Then
\[
\sum_{\mu \in v\Lambda^{\le n_1}} g(s(\mu)) = \sum_{\nu \in
v\Lambda^{\le n_2}} g(s(\nu)).
\]
\end{lemma}
\begin{proof}
Let $n := n_1 \vee n_2$. Using that $v \Lambda^{\le n - n_i}$ is a singleton ($i = 1,2)$ when $v$ lies on an end, one easily checks that $\sum_{\mu \in v\Lambda^{\le n_i}} g(s(\mu)) = \sum_{\lambda \in v\Lambda^{\le n}} g(s(\lambda))$ for $i = 1,2$.
\end{proof}

\begin{proof}[{Proof of Proposition~\ref{prp:sufficient}}]
(a) By definition, graph traces are constant on ends, and hence on
equivalence classes of ends. The formula~\eqref{eq:graph trace
formula} holds by definition of a $k$-graph trace.

(b) Define $\overline{g}(v) : \Lambda^0 \to \R^+$ by
\[
\overline{g}(v) := \sum_{\lambda \in v\Lambda^{\le n_v}} g([x_{s(\lambda)}])
\]
Note that if $x$ is an end of $\Lambda$, then $n_{x(0)} := 0$ has the
property that $s(\lambda)$ lies on an end of $\Lambda$ for each
$\lambda \in x(0)\Lambda^{n_{x(0)}}$. Hence Lemma~\ref{lem:independent of n_v}
shows that $\overline{g}(x(0)) = g([x])$.

Fix $v \in \Lambda^0$ and $n \in \N^k$. We must show that
\begin{equation} \label{eq:must show}
\overline{g}(v) = \sum_{\lambda \in v\Lambda^{\le n}}
\overline{g}(s(\lambda)).
\end{equation}
We may assume without loss of generality that $n \le n_v$ because
if it is not, then $n'_v := n \vee n_v$ can be used in place of
$n_v$ by Lemma~\ref{lem:independent of n_v} and satisfies $n \le
n'_v$. Since $\Lambda^{\le n_v} = \Lambda^{\le n} \Lambda^{\le n_v
- n}$ \cite[Lemma~3.6]{RSY1}, we then have that for each
$\lambda \in v\Lambda^{\le n}$,
the element $n_v - n$ has the property that for each $\alpha \in
s(\lambda) \Lambda^{\le n_v - n}$, the source of $\alpha$ is on an
end of $\Lambda$ and hence
\begin{equation}\label{eq:extension sum}
\overline{g}(s(\lambda)) = \sum_{\alpha \in s(\lambda)\Lambda^{\le
n_v - n}} g(s(\alpha))
\end{equation}
by Lemma~\ref{lem:independent of n_v}. But now
\begin{align*}
\sum_{\lambda \in v\Lambda^{\le n}} \overline{g}(s(\lambda))
&= \sum_{\lambda \in v\Lambda^{\le n}} \Big(\sum_{\alpha \in
s(\lambda)\Lambda^{\le n_v - n}} g(s(\alpha))\Big)
\quad\text{by~\eqref{eq:extension sum}}\\
&= \sum_{\lambda\alpha \in v\Lambda^{\le n_v}}
g(s(\lambda\alpha)) \quad\text{by~\cite[Lemma~3.6]{RSY1}}\\
&= \overline{g}(v)
\end{align*}
by definition of $\overline{g}$.
\end{proof}

Finally we show that we can check that a given function is a graph
trace just by considering edges and vertices in the skeleton of
$\Lambda$. This is useful as it simplifies the task of checking that a
given function is a $k$-graph trace.

\begin{lemma}
Let $\Lambda$ be a locally-convex row-finite $k$-graph. Suppose
that $g : \Lambda^0 \to \R^+$ satisfies $g(v) = \sum_{e \in
v\Lambda^{e_i}} g(s(e))$ for all $v \in \Lambda^0$ and all $1 \le
i \le k$ such that $v\Lambda^{e_i} \not= \emptyset$. Then $g$
is a $k$-graph trace.
\end{lemma}
\begin{proof}
We proceed by induction on $\ell(n) =
\sum_{i=1}^k n_i$. If $\ell(n) = 0$ then $n = 0$ and
$v\Lambda^{\le n} = \{v\}$ for each $v$, so~\eqref{eq:graph
trace} holds trivially. Suppose as an inductive hypothesis that~\eqref{eq:graph trace}
holds for $\ell(n) \le L$, fix $v \in \Lambda^0$
and $n \in \N^k$ with $\ell(n) = L+1$, and write $n = n' + e_i$
where $\ell(n') = L$. By the inductive hypothesis,
\begin{equation}\label{eq:le n-e_i}
g(v) = \sum_{\lambda \in v\Lambda^{\le n'}} g(s(\lambda)),
\end{equation}
and by hypothesis, we know that for each $\lambda \in v\Lambda^{n'}$,
we have
$g(s(\lambda)) = \sum_{e \in s(\lambda)\Lambda^{\le e_i}} g(s(e))$
(if $s(\lambda)\Lambda^{\le e_i} = \{s(\lambda)\}$, this is trivial,
and otherwise it is precisely the hypothesis of the lemma.
Since~\cite[Lemma~3.6]{RSY1} ensures that
$(\lambda,e) \mapsto \lambda e$ is a bijection from
$\{(\lambda,e) : \lambda \in v\Lambda^{\le n'},
e \in s(\lambda)\Lambda^{\le e_i}\}$
to $v\Lambda^{\le n}$, a straightforward calculation shows that
$g(v) = \sum_{\mu \in v\Lambda^{\le n}} g(s(\mu))$.
\end{proof}

\subsection{The $C^*$-algebras of $k$-graphs which admit $k$-graph traces}
In this subsection we give some structural results and $K$-theory
calculations for $C^*(\Lambda)$ when $\Lambda$ is a $k$-graph
which satisfies the hypotheses of
Proposition~\ref{prp:sufficient}.

\begin{prop} \label{prp:bdry path algebra}
Let $\Lambda$ be a $k$-graph, and suppose that the boundary path
$x : \Omega_{k,m} \to \Lambda$ is surjective. Let $v$ denote the
vertex $x(0) \in \Lambda^0$. Then
\begin{itemize}
\item[(a)] the collection $G := \{p - q : p,q \le d(x), x(p) =
x(q)\}$ is a subgroup of $\Z^k$;
\item[(b)] the projection $p_v$
is  full in $C^*(\Lambda)$; and
\item[(c)] there is an
isomorphism $\phi$ of the full corner $p_v C^*(\Lambda) p_v$ onto
the subalgebra $C^*(\{L_n : n \in G\}) \subset C^*(\Z^k)$ which
satisfies $\phi(s_{x(0,p)} s^*_{x(0,q)}) = L_{p-q}$ whenever $x(p)
= x(q)$.
\end{itemize}
In particular, $C^*(\Lambda)$ is Morita equivalent to $C^*(G)$
which is isomorphic to $C(\T^l)$ for some $0 \le l \le k$.
\end{prop}
\begin{proof}
(a) $G$ clearly contains the identity, and is closed under inverses by symmetry. Suppose that $x(p) = x(q)$ and $x(p') = x(q')$, so $n = p-q$ and $n' = p'-q'$ belong to $G$.
We must show that $n + n' \in G$. Since $q, p' \le d(x)$, we have
$q \vee p' \le d(x)$. Let $\alpha := x(q, q\vee p')$ and let
$\beta := x(p', q\vee p')$. Clearly $s(\alpha) = s(\beta)$. But
$r(\alpha) = x(q) = x(p)$, and since $x$ is surjective, it follows
that $x(0, p)\alpha = x(0, p + (q\vee p') - q)$; and similarly, we
have $x(0,q')\beta = x(0, q' + (q \vee p') - p')$. Hence $x(p +
(q\vee p') - q) = x(q' + (q \vee p') - p')$, and so $n + n' = p-q
+ p'-q' = (p + (q\vee p') - q) - (q' + (q \vee p') - p')$ belongs
to $G$ as required.

(b) Since $x$ is surjective, the hereditary subset of
$C^*(\Lambda)$ generated by $v$ as in \cite[\S5]{RSY1} is all of
$\Lambda^0$, so the ideal generated by $p_v$ is $C^*(\Lambda)$ as
required.

(c) We have that $p_v C^*(\Lambda) p_v = \clsp\{s_\lambda s^*_\mu : \lambda,\mu \in v\Lambda, s(\lambda) = s(\mu)\} = \clsp\{s_{x(0,p)} s^*_{x(0,q)} : x(p) = x(q)\}$ because $x$ is surjective. Suppose that $x(p) = x(q)$ and $x(p') = x(q')$, and that $p-q = p'-q'$. Using~(CK4) and that $x$ is surjective, one checks that $s_{x(0,p)} s^*_{x(0,q)} = s_{x(0,p')} s^*_{x(0,q')}$, so for $n \in G$ we may define $U_n := s_{x(0,p)} s^*_{x(0,q)}$ for any $p,q$ such that $x(p) = x(q)$ and $p-q = n$.

The Cuntz-Krieger relations show that $U_n U^*_n = U^*_n U_n = p_v$ for all $n$, so the $U_n$ are unitaries.

The calculation
\[
U_{q - p} = s_{x(0,q)} s^*_{x(0,p)} = (s_{x(0,p)}
s^*_{x(0,q)})^* = U^*_{p-q}.
\]
Shows that $U_{-n} = U^*_n$ for all $n \in G$. Moreover, since $x$ is surjective, we have $s^*_{x(0,q)} s_{x(0,p')} = s_{x(q, q \vee p')} s^*_{x(p', q\vee p')}$ by \cite[Proposition~3.5]{RSY1} and~(CK3), and it follows that
\[
U_{p-q} U_{p'-q'}
= s_{x(0,p)} s^*_{x(0,q)} s_{x(0,p')} s^*_{x(0,q')}
= s_{x(0, p + (q \vee p') - q)} s^*_{x(0, q' + (q \vee p') - p')}
= U_{p + p' - q - q'}
\]
so that $n \mapsto U_n$ is a representation of $G$. It follows that there is a surjective $C^*$-homomorphism $\phi : C^*(G) \to C^*(\{U_n : n \in G\}) = p_v C^*(\Lambda) p_v$ which satisfies $\phi(\chi_n) = U_n$ for all $n \in G$.

It remains only to show that $\phi$ is injective. For this, we need only
show that $U_m \not= U_n$ for $m \not= n$ and that each $U_n$ where $n \not= 0$ has
full spectrum. Since $\gamma_z(U_m U^*_n) = z^{m-n}U_m U^*_n$, an appropriate choice of $z \in \TT^k$ shows that $U_m \not= U_n$ for $m \not= n$. That each $U_n$ has full spectrum follows from an argument identical to that used in \cite[Lemma~3.9]{PRRS}. This establishes~(c).

The Morita equivalence of $C^*(\Lambda)$ with $C^*(G)$ follows immediately from (2)~and~(3), and since $G$ is a subgroup of $\Z^k$, we must have $G \cong \Z^l$ for some $0 \le l \le k$.
\end{proof}

For the remainder of the appendix, $\Lambda$ will be a fixed
$k$-graph which satisfies the hypotheses of
Proposition~\ref{prp:sufficient}. \vskip.25em

A vertex $v \in \Lambda^0$ lies on an end of $\Lambda$ if and only if $|v\Lambda^{\le n}| = 1$ for all $n \in \N^k$. Let $\Ends(\Lambda)^0$ denote the collection of all such vertices; for each $v \in \Ends(\Lambda)^0$, there is a unique end $x(v)$ whose range is $v$. By a simple argument, we may select a set $V \subset \Ends(\Lambda)^0$ such that for each $x \in \Ends(\Lambda)$ there is a unique $v \in V$ such that $x \sim x(v)$.  We fix this collection for the remainder of the section.

\begin{prop}\label{prp:direct sum}
For each $v \in V$, let  $\Lambda(v)$ be the image of $x(v)$ which is a
subcategory of $\Lambda$. Then each $(\Lambda(v), d|_{\Lambda(v)})$ is itself a $k$-graph, and $C^*(\Lambda)$ is Morita equivalent to $\bigoplus_{v \in V} C^*(\Lambda(v))$.
\end{prop}
\begin{proof}
Since each $x(v)$ is an end, each $\Lambda(v)^0$ is a hereditary subset of $\Lambda^0$. For distinct $v,w \in V$, we have $\Lambda(v) \cap \Lambda(w) = \emptyset$ because otherwise $x(v) \sim x(w)$ contradicting our choice of $V$. By \cite[Theorem~5.2]{RSY1},
for each $v \in \Lambda^0$, the projection $P_v := \sum_{w \in \Lambda(v)^0} p_w$ determines an ideal $I_v := C^*(\Lambda) P_v C^*(\Lambda)$ which is Morita equivalent to $C^*(\Lambda(v)^0 \Lambda) = C^*(\Lambda(v))$. Each $I_v = \clsp\{s_\lambda s^*_\mu : s(\lambda) = s(\mu) \in \Lambda(v)^0\}$, and since the distinct $\Lambda(v)$ do not intersect the ideal generated by all the $P_v$ is isomorphic to $\bigoplus_{v \in V}
I_v$.

The assumption that for each vertex $w \in \Lambda^0$ there is an element $n_w \in \N^k$ such that $s(w\Lambda^{\le n_w}) \subset \Ends(\Lambda)^0$ guarantees that every vertex of $\Lambda$ belongs to the saturated hereditary set generated by $V$. Another application of \cite[Theorem~5.2]{RSY1} shows that the ideal generated by all the $I_v$ is $C^*(\Lambda)$. Hence $C^*(\Lambda) = \bigoplus I_v$ is Morita equivalent to
$\bigoplus C^*(\Lambda(v))$.
\end{proof}

\begin{cor} \label{cor:k-theory for ends}
For each $v \in V$, let $G_v := \{p - q : x(v)(p) = x(v)(q)\} \subset \Z^k$. Then $C^*(\Lambda)$ is Morita equivalent to $\bigoplus_{v \in V} C^*(G_v) \cong \bigoplus_{v \in V} C(\T^{l_v})$ where $0 \le l_v \le k$ for each $v$. In particular $K_*(C^*(\Lambda))$ is isomorphic to $\bigoplus_{v \in V} K_*(C(\T^{l_v}))$.
\end{cor}

\end{document}